\documentclass[sn-mathphys]{sn-jnl}
\usepackage{graphicx}%
\usepackage{longtable,tabularx}
\usepackage{amsmath,amssymb,amsfonts}%
\usepackage{amsthm}%
\usepackage[referable]{threeparttablex}
\usepackage{mathrsfs}%
\usepackage[title]{appendix}%
\usepackage{xcolor}%
\usepackage{textcomp}%
\usepackage{adjustbox}
\usepackage{manyfoot}%
\usepackage{booktabs}%
\usepackage{amsmath}
\usepackage{algorithm}%
\usepackage{algorithmicx}%
\usepackage{algpseudocode}%
\usepackage{arydshln}
\usepackage{listings}%
\usepackage{nomencl}
\makenomenclature
\usepackage{subfig}
\usepackage{graphicx}
\usepackage{amsmath}
\usepackage[version=4]{mhchem}
\usepackage{siunitx}
\usepackage{longtable,tabularx}
\usepackage{hyperref}
\usepackage{dsfont}
\def\abs#1{\left\lvert#1\right\rvert}
\newcommand{\eqnref}[1]{Eq.~(\ref{#1})}
\newcommand{\figref}[1]{Fig.~\ref{#1}}
\newcommand{\tabref}[1]{Table~\ref{#1}}
\usepackage [autostyle,english = american]{csquotes}
\MakeOuterQuote{"}
\usepackage{tikz}
\usepackage{forest}
\usepackage{stanli}
\usetikzlibrary{positioning,shapes,shadows,arrows}
\newtheorem{theorem}{Theorem}

\newtheorem{example}{Example}

\newcommand{\inp}{\text{in}}
\newcommand{\outp}{\text{out}}

\newcommand{\tstar}[5]{% inner radius, outer radius, tips, rot angle, options
\pgfmathsetmacro{\starangle}{360/#3}
\draw[#5] (#4:#1)
\foreach \x in {1,...,#3}
{ -- (#4+\x*\starangle-\starangle/2:#2) -- (#4+\x*\starangle:#1)
}
-- cycle;
}
\begin{document}

\newcommand{\Title}%{Mixed categorical Gaussian process and Bayesian optimization for high-dimensional multidisciplinary applications}
{High-dimensional mixed-categorical Gaussian processes with application to multidisciplinary design optimization for a green aircraft}

\title[\Title]{\Title}

\author*[1,2]{\fnm{Paul} \sur{Saves}}\email{paul.saves@onera.fr}
\author[3]{\fnm{Youssef} \sur{Diouane}}\email{youssef.diouane@polymtl.ca}
\author[1]{\fnm{Nathalie} \sur{Bartoli}}\email{nathalie.bartoli@onera.fr}
\author[1]{\fnm{Thierry} \sur{Lefebvre}}\email{thierry.lefebvre@onera.fr}
\author[4]{\fnm{Joseph} \sur{Morlier}}\email{joseph.morlier@isae-supaero.fr}

\affil*[1]{ONERA/DTIS, Université de Toulouse, Toulouse, France}
\affil[2]{ISAE-SUPAERO, Université de Toulouse, Toulouse, France}
\affil[3]{GERAD and Department of Mathematics and Industrial Engineering, Polytechnique Montr\'eal, Montréal, QC, Canada}
\affil[4]{ICA, Universit\'e de Toulouse, ISAE--SUPAERO, INSA, CNRS, MINES ALBI, UPS, Toulouse, France}

\abstract{
%Multidisciplinary design optimization methods aim to adapt numerical optimization techniques to the design of engineering systems involving multiple disciplines. In this context, a large number of mixed continuous, integer, and categorical variables might arise during the optimization process, and practical applications involve a significant number of design variables.
\textcolor{black}{
Recently, there has been a growing interest in mixed-categorical metamodels based on Gaussian Process (GP) for Bayesian optimization. In this context,  different approaches can be used to build the mixed-categorical GP. Many of these approaches involve a high number of hyperparameters; in fact, the more general and precise the strategy used to build the GP, the greater the number of hyperparameters to estimate. This paper introduces an innovative dimension reduction algorithm that relies on partial least squares regression to reduce the number of hyperparameters used to build a mixed-variable GP. Our goal is to generalize classical dimension reduction techniques commonly used within GP (for continuous inputs) to handle mixed-categorical inputs. The good potential of the proposed method is demonstrated in both structural and multidisciplinary application contexts. The targeted applications include the analysis of a cantilever beam as well as the optimization of a green aircraft, resulting in a significant 439-kilogram reduction in fuel consumption during a single mission.
}}

\keywords{Surrogate modeling, Gaussian process, Mixed-categorical inputs, High-dimension, {\color{black}Bayesian optimization}, Aircraft design}
\maketitle

%\tableofcontents
%\newpage
\newpage
\color{black}
\section*{Nomenclature}

\begin{tabular}{lcl}
GP & : &Gaussian Process \\
DoE & : &Design of Experiments \\  
LHS & : &Latin Hypercube Sampling \\
PLS & : &Partial Least Squares \\
KPLS & : &Kriging with Partial Least Squares \\
BO & : &Bayesian Optimization \\
%SBO & : &Surrogate-Based Optimization \\
$n$ & : &number of continuous variables \\
$m$ & : &number of integer variables \\
$l$ & : &number of categorical variables \\
$\Omega \in \mathbb{R}^n$ &: & continuous space\\
$S \in \mathbb{Z}^m$  &: & integer space \\
$\mathbb{F}^l$ & : &categorical space \\
$L_i $ & : &number of levels for the  $i^{th}$ categorical variable \\
$\theta^{cont}$  &  : & hyperparameters related to the continuous variables for GP  \\
$\theta^{int}$  &  : &hyperparameters related to the integer variables for GP  \\
%$k$ & correlation kernel \\
%$\theta^{cont}_j, j \in \{1,\ldots, n+m\} $  &  hyperparameter  for the  $j^{th}$ continuous or integer variable  \\ 
$\theta^{cat}$  &: &hyperparameters related to the categorical variables for GP  \\
$\Theta_i , $  &  : & matrix of hyperparameters  for the  $i^{th}$ categorical variable  \\
$R^{cont} (\theta^{cont})$ & : &correlation matrix for continuous inputs \\
$R^{int} (\theta^{int}) $ & : &correlation matrix for integer inputs \\
$R^{cat} (\theta^{cat}) $ & : &correlation matrix for categorical inputs \\ 

$\Theta=\{ \theta^{cont},\theta^{int},\theta^{cat}\}$ &  : &hyperparameters for the Gaussian process model  \\ 
$G_*$ & : &PLS rotation matrix \\
$\zeta_{c^r_i}$ & : &full space cross-level relaxation of a categorical variable \\
$\hat{\theta}$ & : &continuous or integer hyperparameters for PLS approximation \\
$\hat{\Theta}$ & : &categorical hyperparameters for PLS approximation \\
MDA & : &Multidisciplinary Design Analysis \\
MDO & : &Multidisciplinary Design Optimization \\
EGO & : &Efficient Global Optimization \\
GD & : &Gower Distance \\
CR & : &Continuous Relaxation \\
HH & : &Homoscedastic Hypersphere \\
EHH & : &Exponential Homoscedastic Hypersphere \\
SPD & : &Symmetric Positive Definite \\
MLE & : &Maximum Likelihood Estimator \\
RMSE & : &Root Mean Squares Error \\
SMT & : &Surrogate Modeling Toolbox \\
SEGOMOE & : &Super Efficient Global Optimization with Mixture Of Experts \\
FAST-OAD & : &Future Aircraft Sizing Tool with Overall Aircraft Design \\
MAC & : &Mean Average Chord \\
TOFL & : &Take-Off Field Length \\
VT & : &Vertical Tail \\
HT & : &Horizontal Tail \\

DRAGON & : & Distributed fans Research Aircraft with electric Generators  \\
     & &   by ONera  \\
\end{tabular}
\color{black}

\section{Introduction}
\label{SMO_sec:intro}

Costly black-box simulations play an important role for many engineering and industrial applications. 
For this reason, surrogate modeling has been extensively used across a wide range of use cases, including aircraft design~\cite{SciTech_cat}, deep neural networks~\cite{snoek2015scalable}, coastal flooding prediction~\cite{lopez}, agricultural forecasting~\cite{MLP}, and seismic imaging~\cite{YDiouane_SGratton_XVasseur_LNVicente_HCalandra_2016}.
These black-box simulations are generally complex and may involve mixed-categorical input variables. For instance, a Multidisciplinary Design Analysis (MDA) aircraft design tool~\cite{David_2021} must consider mixed variables such as the number of engines or the list of possible materials~\cite{SciTech_cat}.

In this paper, our objective is to develop an affordable surrogate model, denoted as~$ \hat{f}$, for a black-box function that involves mixed variables given by
\begin{equation}
f :  \Omega \times S \times \mathbb{F}^l \to \mathbb{R}.
  \label{SMO_eq:opt_prob}
\end{equation}
This function $f$ represents a computationally expensive black-box simulation.
$\Omega \subset \mathbb{R}^n$ denotes the bounded continuous design set for the $n$ continuous variables.  $S \subset \mathbb{Z}^m$ denotes the bounded integer set where $L_1, \ldots, L_m$ are the numbers of levels of the $m$ quantitative integer variables on which we can define an order relation and \hbox{$ \mathbb{F}^l = \{1, \ldots, L_1\} \times \{1, \ldots, L_2\} \times  \ldots \times \{1, \ldots, L_l\}$} is the design space for the $l$ categorical qualitative variables with their respective  $L_1, \ldots, L_l$ levels.

For such purpose, the use of a Gaussian Process (GP)~\cite{williams2006gaussian}, also called Kriging model~\cite{krige1951statistical}, is recognized as an effective modeling approach for constructing a response surface model based on an available dataset.  Specifically, we make the assumption that our unknown black-box function, denoted as $f$, follows a \textcolor{black}{GP} with mean $\mu^{f}$ and standard deviation $\sigma^f$, expressed as follows:
\begin{equation}
f \sim \hat{f}=
%\mbox{GP}
\mathcal{GP}
\left(\mu^{f}, (\sigma^f)^2\right). 
\label{SMO_eq:GP:f}
\end{equation}
Several modeling approaches have been put forward for addressing the challenges of handling categorical or integer variables within the context of \textcolor{black}{GP}~\cite{Pelamatti, Zhou, Deng, Roustant, GMHL, Gower, cuesta2021comparison, SciTech_cat}. In comparison to \textcolor{black}{GP} designed for continuous variables, the most important changes concern the estimation of the correlation matrix, an essential element in the derivation of $\mu^{f}$ and $\sigma^f$. Much like the procedure for constructing a \textcolor{black}{GP} with continuous inputs, Continuous Relaxation (CR) techniques~\cite{GMHL, SciTech_cat}, models involving continuous latent variables~\cite{cuesta2021comparison}, and Gower Distance (GD) based models~\cite{Gower} use a kernel-based approach for estimating this correlation matrix.
However, recent innovative approaches take a different path by modeling directly the various entries of the correlation matrix~\cite{Pelamatti,Zhou,Deng,Roustant}, and therefore, do not rely on any kernel choice\textcolor{black}{, such methods involve the Homoscedastic Hypersphere (HH)~\cite{Zhou} and the Exponential HH (EHH)~\cite{Mixed_Paul} kernels.} It has been shown in~\cite{Mixed_Paul} that the \textcolor{black}{HH} correlation kernel generalizes simpler methods like CR or GD kernels. However, this more general method for handling categorical design variables increases the number of hyperparameters required to be tuned associated with the GP model. In particular, this means that the method could only be used for small dimensional problems. 

\textcolor{black}{Many efficient approaches have been proposed for handling a high number of continuous variables within GP~\cite{Bouhlel18, bouhlel_KPLSK, bouhlel2019gradient}. The Kriging with Partial Least Squares (KPLS) method~\cite{Bouhlel18,bouhlel_KPLSK} is one of the most commonly used reduction techniques~\cite{KPLSu1,KPLSu2} to tackle high dimensional data. Several other dimension reduction methods include principal components analysis~\cite{wang2017batched}, polynomial chaos expansion~\cite{zuhal2021dimensionality}, radial basis functions~\cite{regis2020survey}, active subspace~\cite{AS}, manifold embedding~\cite{tenenbaum2000global} or marginal Gaussian process~\cite{MGP}. The KPLS technique allows constructing the GP model with the same continuous variables but using a few number of parameters; which reduces significantly the computational cost of computing a GP model. }
\textcolor{black}{
For mixed-categorical GP models, given that the computational effort related to the construction of the GP model may not scale well to practical applications involving categorical variables with a large number of levels, the number of hyperparameters to be tuned need to be considered more thoroughly. 
In the literature, GPs have been applied to no more than 15 hyperparameters due to the high computational cost associated with the estimation of the hyperparameters~\cite{GP14}. Adapting dimension reduction techniques, such as KPLS, to mixed-categorical GPs will thus enable solving practical mixed-categorical engineering problems where often a high number of hyperparameters is required to be estimated. To the best of our knowledge, there is no equivalent approach to handle mixed-categorical data without using relaxation techniques. All existing dimension reduction techniques, including KPLS, are not adapted for advanced mixed-categorical GP models such as HH or EHH. 
We note also that, although this paper focuses mainly on surrogate modeling, the proposed models can be integrated within any surrogate-based optimization method~\cite{AuHa2017}, such as surrogate-based evolutionary algorithm~\cite{sbea1,sbea2} or a Bayesian Optimization (BO) method~\cite{Jones2001JOGO}.}  

%\textcolor{black}{Concerning the latent variables approach~\cite{zhang2020latent}, it has not been investigated because of three major limitations. First, this method assumes a continuous order relation between the levels of the categorical variables, which is something we seek to be rid of to be as general as possible. Second, such an order relation and use of continuous kernel imply positive correlation, which is not always the case for categorical variables. Third, such a method embeds the categorical correlations in smaller subspaces but keeps the number of hyperparameters associated with the GP model high ( $l$ different 2D continuous latent subspaces but $ \Sigma_{i=1}^l \left( 2 L_i -3 \right)$ parameters to be tuned). More recently, the latent variable approach have been generalized and the first limitation has been removed thanks to a one-hot encoding~\cite{oune2021latent}. Another improvement brought by~\cite{oune2021latent} is that it requires a unique 2D latent space but this model still requires $ \Sigma_{i=1}^l \left(  2 L_i \right)$ parameters to be tuned. }

\textcolor{black}{
In this work, we target to use dimension reduction techniques for reducing the number of hyperparameters within the GP in order to allow modeling efficiently high-dimension mixed-categorical data. In this context, high dimensionality is related to the high number of categorical variables potentially with a high number of levels (a few dozen). In fact, using relaxation approaches (by converting categorical choices to continuous variables) leads to a very high number of hyperparameters to estimate, particularly for high resolution approaches such as those based on HH and EHH kernels.
%In this work, we show how to reduce the computational cost related to the construction of the mixed categorical GP model proposed in~\cite{Pelamatti} using KPLS. 
We have also specifically used our proposed mixed-categorical GP models, within a BO framework, to solve a constrained optimization problem involving expensive-to-compute black-box simulations for objective and constraints functions~\cite{Martins2021}.
The proposed approach is shown in particular to be efficient in solving a high dimensional mixed-categorical Multidisciplinary Design Optimization (MDO) problem~\cite{Lambe2012}. }
\textcolor{black}{All the GP models proposed in this work are implemented in the open-source Surrogate Modeling Toolbox (SMT)\footnote{\url{https://smt.readthedocs.io/en/latest/}}~\cite{saves2023smt}. }
%More precisely, this work extends HH, the most general categorical kernel formulation, to high-dimension using partial least squares technique. This new high-dimensional mixed-categorical kernel generalizes existing high-dimensional continuous kernels where KPLS is often coupled with CR to handle categorical variables. Additionally, our proposed method allows us to recover an approximate estimation of the correlation matrix between the levels of a given categorical variable, allowing for further investigation and analysis.}

%This toolbox has been used regularly in engineering problems, for aircraft engine consumption modeling~\cite{DGP1}, multi-fidelity sensitivity analysis~\cite{drouet2023multi}, high-order robust finite elements methods~\cite{karban2021fem,kudela2022recent}, chemical process design~\cite{savage2020adaptive} and many other applications. These applications are generally multidisciplinary oriented because surrogate models are of high interest in such settings.}

The remainder of this paper is the following.
In Section~\ref{SMO_sec:GP}, a detailed review of the GP model for continuous and for categorical inputs is given. In Section~\ref{SMO_sec:PLS_KPLS}, we present the PLS regression for vectors and matrices and their application to GP model for both continuous and categorical inputs. 
Section~\ref{SMO_sec:Results} presents academical tests as well as the obtained results on multidisciplinary optimization.
Conclusions and perspectives are finally drawn in Section~\ref{SMO_sec:conclu}.

\noindent\textit{Notations:} For a vector $x$, both notations $[x]_j$ and $x_j$ stand  for the $j^{th}$ component of $x$. Similarly, the $i$ (row index) and $j$ (column index) entry of a matrix $X$ is denoted $[X]^{j}_{i}$.  
\section{GP for mixed-categorical inputs}
\renewcommand{\footnoterule}{ \hrule width5cm \vspace*{0.1cm} }    
\label{SMO_sec:GP}

In this section, we present the mathematical background associated with GP for mixed-categorical variables. This part also introduces the notations  used throughout the paper. Here, the general case involving mixed integer variables is considered. Namely, we assume that $f:\mathbb{R}^n \times  \mathbb{Z}^m \times \mathbb{F}^l \mapsto \mathbb{R}$ and our goal is to build a GP surrogate model for $f$. 

\subsection{A mixed GP formulation} 
Given a set of data points, called a Design of Experiments (DoE)~\cite{forrester}, Bayesian inference learns the GP model that explains the best this dataset. A GP model consists of a mean response hypersurface $\mu^{f}$, as well as an estimation of its variance $(\sigma^f)^2$. In the following, $n_t$ denotes the size of the given DoE dataset $(W, \textbf{y}^f)$ such that $W=\{w^1,w^2,\ldots,w^{n_t}\} \in (\mathbb{R}^n \times  \mathbb{Z}^m \times \mathbb{F}^l)^{n_t}$ and $\textbf{y}^f=[f(w^1),f(w^2),\ldots,f(w^{n_t})]^{\top}$. 
For an arbitrary $w= (x,z,c) \in \mathbb{R}^n \times  \mathbb{Z}^m \times \mathbb{F}^l$, not necessary in the DoE, the GP model prediction at $w$ writes as $\hat f(w) = \mu ({w})+\eta(w) \in \mathbb{R} $, with $\eta$ being the uncertainty between $\hat{f}$ and the model approximation $\mu$~\cite{GP14}. The considered error terms are random variables of variance $\sigma^2$.  Using the DoE, the expression of  $\mu^{f}$ and the estimation of its variance $(\sigma^f)^2$ are given as follows:

\begin{equation} \label{SMO_eq:mean:GP}
\mu^f(w)= \hat{\mu}^f+r(w)^\top  [R(\Theta)]^{-1}({y}^f-\mathds{1} \hat{\mu}^f), 
\end{equation}
and
\begin{equation}
\label{SMO_eq:std:GP}
(\sigma^f(w))^2=(\hat{\sigma}^f)^2\left(1-r(w)^\top  [R(\Theta)]^{-1}r(w)+ \frac{ \left(1-\mathds{1}^\top  [R(\Theta)]^{-1}r(w) \right)^2}{\mathds{1}^\top  [R(\Theta)]^{-1}\mathds{1}}\right), \end{equation}
where $\hat{\mu}^f$ and $\hat{\sigma}^f$ are, respectively, the Maximum Likelihood Estimator (MLE)~\cite{MLE} of $\mu$ and $\sigma$. $\mathds{1}$ denotes the vector of $n_t$ ones. $R$ is the $ n_t \times n_t $ correlation matrix between the input points and $r(w)$ is the correlation vector between the input points and a given $w$.
To have a compact notation, let $[A]_i^j$ denote the coefficient of the matrix $A$ in the $i^{\text{th}}$ row and $j^{\text{th}}$ column.
The correlation matrix $R$ is defined, for a given couple of indices  $(r,s) \in \{1,\ldots,n_t\}^2$, by \begin{equation}
\label{SMO_eq:R}
 [R(\Theta)]_{r}^{s}=k\left(w^r,w^s,\Theta\right) \in \mathbb{R},\end{equation}
and the vector $r(w)\in \mathbb{R}^{n_t}$ is defined as $ r(w) =[k(w,w^1,\Theta), \ldots , k(w,w^{n_t},\Theta)]^{\top}$,
where $k$ is a given correlation kernel that relies on a set of hyperparameters $\Theta$~\cite{Roustant}. 
The mixed-categorical  correlation kernel is given as the product of three kernels:
\begin{equation}
k(w^r,w^s,\Theta) =  k^{cont}\left(x^r,x^s,\theta^{cont}\right) k^{int}\left(z^r,z^s,\theta^{int}\right)
k^{cat}\left(c^r,c^s,\theta^{cat}\right),
\label{SMO_eq:decomp_mix}
\end{equation}
where $k^{cont}$  and $\theta^{cont}$ are respectively the continuous kernel and its associated hyperparameters, $k^{int}$  and $\theta^{int}$ are the integer kernel and its hyperparameters, and last $k^{cat}$  and $\theta^{cat}$ are the ones related with the categorical inputs. In this case, one has $\Theta=\{ \theta^{cont},\theta^{int},\theta^{cat}\}$. 
Henceforth, the general correlation matrix $R$ will rely only on the set of the hyperparameters $\Theta$:
\begin{equation}
    \label{SMO_eq:corel:mat}
    [R(\Theta)]_{r}^{s} = [R^{cont}(\theta^{cont})]_{r}^{s}  [R^{int}(\theta^{int})]_{r}^{s}
    [R^{cat}(\theta^{cat})]_{r}^{s},
\end{equation}
where $[R^{cont}(\theta^{cont})]_{r}^{s} =k^{cont}(x^r,x^s,\theta^{cont}) $, $[R^{int}(\theta^{int})]_{r}^{s} =k^{int}(z^r,z^s,\theta^{int}) $ and  $[R^{cat}(\theta^{cat})]_{r}^{s}=k^{cat}(c^r,c^s,\theta^{cat})$. 
The set of hyperparameters $\Theta$ could be estimated using the DoE dataset $({W},{y}^f)$ through the MLE approach on the following way

\begin{equation}
\Theta^*= \arg\max_{\Theta} \mathcal{L}(\Theta)=\left( - \frac{1}{2} {{y}^f}^\top [R(\Theta)]^{-1} {{y}^f}   - \frac{1}{2} \log 	\abs{  [R(\Theta)]} - \frac{n_t}{2} \log 2 \pi    \right),
\label{SMO_eq:likelihood}
\end{equation}
where $R(\Theta)$ is computed using~\eqnref{SMO_eq:corel:mat}.
To construct the correlation matrix for continuous or integer inputs, several choices for the correlation kernel are possible. Usual families of kernels include exponential kernels or Matern kernels~\cite{Lee2011}. 
In contrast, to construct the correlation matrix for categorical inputs, we can either use a kernel as for the continuous or integer variables or we can directly model the entries of the correlation matrix thanks to a Symmetric Positive Definite (SPD) matrix parameterization. The latter approach is what is done for the HH kernel, for example~\cite{Roustant}. For this kernel, the hyperparameters $\theta^{cat}$ can be seen as a concatenation of a set of symmetric matrices, \textit{i.e.}, $\theta^{cat} = \{  \Theta_1, \Theta_2, \ldots, \Theta_l \} $. The construction of this kernel is thus relying on the estimation of $\sum_{i=1}^l \frac{1}{2} L_i (L_i-1) $ hyperparameters.

\subsection{The homogeneous categorical kernel} 
A recent paper~\cite{Mixed_Paul} unified the kernel-based approach and the matrix-based approach through the homogeneous model described hereafter.
Recall that $l$ denotes the number of categorical variables. 
For a given $i \in \{1, \ldots, l\}$, let $c^r_{i} $  and $c^s_{i} $  be a couple of categorical variables taking respectively the $\ell^r_i$ and the $\ell^s_i$ level on the categorical variable $c_i$.
The hyperparameter matrix peculiar to this variable $c_i$ is 
$$\Theta_i= \begin{bmatrix}
[\Theta_i]_{1}^{1} & \textcolor{white}{9} & \hspace{2em} { \textbf{\textit{ Sym.}}}  \textcolor{white}{9} & \\
[\Theta_i]_{1}^{2}  & [\Theta_i]_{2}^{2} & \textcolor{white}{9} \\
\vdots &\ddots & \ddots & \textcolor{white}{9}  \\
[\Theta_i]_{1}^{L_i} &  \ldots & [\Theta_i]_{L_i-1}^{L_i} &[\Theta_i]_{L_i}^{L_i} \\ 
\end{bmatrix}.$$
First, the correlation term $[R^{cat}(\theta^{cat})]_{r}^{s}$ can be formulated in a level-wise form~\cite{Pelamatti} as: 
\begin{equation}
\begin{split}
k^{cat}(c^r,c^s,\theta^{cat}) =  &{\displaystyle \prod_{i=1}^{l}  [R_i(\Theta_i)]_{\ell^r_i}^{\ell^s_i} } \\
= &\ {\displaystyle \prod_{i=1}^{l}  \kappa ( 2 [ \Phi(\Theta_i) ]_{{ \ell_i^r}}^{{\ell_i^s}} ) \  \kappa ( [ \Phi(\Theta_i) ]_{{ \ell_i^r}}^{{\ell_i^r}} ) \  \kappa ( [ \Phi(\Theta_i) ]_{{ \ell_i^s}}^{{\ell_i^s}} ), }
\end{split}
\label{SMO_eq:homo_HS}
\end{equation}
where $\kappa(.)$ is either a positive definite kernel or the identity and $\Phi(.)$ is a SPD function such that the matrix $\Phi(\Theta_i) $ is SPD if $\Theta_i$ is SPD.
For an exponential kernel,~\tabref{SMO_tab:kernels} gives the parameterizations of $\Phi$ and $\kappa$ that correspond to GD, CR, HH and EHH kernels. 
Note that the complexity between these different kernels is reflected by the number of hyperparameters that characterize them.
As in~\cite{Mixed_Paul}, for all categorical variables $i \in \{1, \ldots, l\}$, the matrix $C(\Theta_i)\in \mathbb{R}^{L_i \times L_i}$ is lower triangular and built using a hypersphere decomposition~\cite{HS,HS_Jacobi} from the symmetric matrix $\Theta_i \in \mathbb{R}^{L_i \times L_i}$ of hyperparameters. 
The variable $\epsilon$ is a small positive constant and the variable $\theta_{i}$ denotes the only positive hyperparameter that is used for the GD kernel.
Nevertheless, until now, PLS regression was only applied to mixed integer inputs for the CR kernel~\cite{SciTech_cat}. In the following section, we will show how to extend the PLS regression for the more general HH kernel. 
%
%\resizebox{0.97\textwidth}{!}{%
    \begin{table}[th]
    \caption{Kernels based on~\eqnref{SMO_eq:homo_HS} formulation.}
    \centering
    \begin{tabular*}{\linewidth}{ccll}
    \hline
    \textbf{Name} & $\kappa(\phi)$   &  \hspace{3cm} ${\centering \Phi(\Theta_i)}$   &  \hspace{-0.5cm}  \# of Hyperparam. \\
    \hline
    \textbf{GD}   &  $\exp(-\phi) $ & ${ \scriptstyle [\Phi(\Theta_i)]_{j}^{j} :=  \frac{1}{2}  \theta_{i} \quad ~;~ [\Phi(\Theta_i)]_{j \neq j'}^{j'} := 0 }$ & 1   \\
     \textbf{CR}  & $\exp(-\phi) $ &  $  { \scriptstyle [\Phi(\Theta_i)]_{j}^{j} := [\Theta_i]_{j}^{j} \  ~;~ [\Phi(\Theta_i)]_{j \neq j'}^{j'} := 0 } $  & $L_i$  \\
     \textbf{EHH}  & $\exp(-\phi)$ & 
    $  { \scriptstyle [\Phi(\Theta_i)]_{j}^{j} := 0 \quad \quad ~;~ [\Phi(\Theta_i)]^{j'}_{j \neq j'} := \frac{\log \epsilon }{2} ([C(\Theta_i) C(\Theta_i) ^\top]_{j}^{j'} -1)  }$  & $\frac{1}{2}  (L_i)  (L_i-1) $\\
     \textbf{HH}  &  $\phi$ &    $  { \scriptstyle [\Phi(\Theta_i)]_{j}^{j} := 1 \quad \quad ~;~ [\Phi(\Theta_i)]_{j \neq j'}^{j'} := \frac{1}{2} [C(\Theta_i) C(\Theta_i)^\top]_{j}^{j'} }$ & $\frac{1}{2}  (L_i)  (L_i-1) $  \\
      \hline
    \end{tabular*}%
\label{SMO_tab:kernels}
\end{table}
%} 
%\vspace{0.3cm}

\section{KPLS for mixed-categorical inputs}
\label{SMO_sec:PLS_KPLS}

To have a sparse model that can extend to high dimension, and to facilitate the optimization of the hyperparameters, one seeks to express the correlation matrix  $R^{cont} (\theta^{cont})$ with $d \ll n $ relevant hyperparameters. Such a method is KPLS~\cite{Bouhlel18} that is an adaptation of the Partial Least Squares (PLS) regression for exponential kernels. To introduce the variables and notations, the next part presents a short recall of PLS regression for vector inputs and of KPLS for continuous variables. Then, the second part presents our extension to matrix inputs and its application for categorical variables.

\subsection{KPLS for continuous inputs}

In this part, we introduce the PLS regression for vector inputs as it was developed by Wold~\cite{wold_1975} and its application to GP kernels (namely KPLS) developed by Bouhlel~\textit{et al.}~\cite{Bouhlel18}.

\subsubsection{PLS for vectors inputs}
\label{SMO_subsec:PLS_cont}

We present here the classical method for the continuous case but integer inputs can be treated similarly by considering them as continuous. Let the DoE be $({X},{y}^f)$ where $X$ is the continuous data matrix of size $n_t \times n$ and ${y}^f$ is the response vector of size $n_t$. The PLS regression method is designed to search out the best multidimensional direction in $\mathbb{R}^n$ that explains the output ${y^f}$~\cite{wold_1975}.
The first principal component (or \textit{score}) $h^{(1)}$ is obtained by searching the best direction  (or \textit{weight}) $g^{(1)}$ that maximizes the squared covariance between $h^{(1)} =  X  g^{(1)} $  and ${y^f}$: 
\begin{equation}
    g^{(1)}= \arg \max_{g \in \mathbb{R}^n} \left\{ g^\top X^\top y^f {y^f}^\top X g ~~\mbox{s.t.}~~ g^\top g =1  \right\}.
    \label{SMO_eq:catKPLS}
\end{equation}
Next, we compute the residuals of the inputs as $X^{(1)}= X - \xi^{(1)} h^{(1)}   $ with $\xi^{(1)}$ being the regression coefficients (or \textit{loadings}) that minimize the residual for every point.  We also project the output and, denoting $\gamma_1$ the corresponding coefficient, we have $y^{f(1)}= y^f - \gamma_1 h^{(1)}  $.
For all $t \in \{1, \ldots, d \} $, we iterate the process of~\eqnref{SMO_eq:catKPLS} with the residuals $X^{(t)}$ and $y^{f(t)}$. 
At the end of the process, we can write the various computed quantities in a matrix form. Namely, we denote $G$ the $n \times d$ matrix such that $g^{(t)}$ is the $t^{th}$ column of $G$ and $\Xi$ the $n \times d$ matrix such that $\xi^{(t)}$ is the $t^{th}$ column of $\Xi$.
Let $G_* = G(\Xi^T G)^{-1}$ be the $n \times d$ matrix such that $G_* = [G_*^{1}, \ldots, G_*^{d}]$ with $G_*^{t} \in \mathbb{R}^n    $, $G_*$ is called the rotation matrix~\cite{Bouhlel18,PLS-based}. Thanks to this matrix, we can express the score $h^{(t)}$ as a function of the input $X$ as: $$ \forall t \in \{1, \ldots, d \}, h^{(t)}=X^{(t)} g^{(t)}=X G_*^{t}.$$ 
By PLS, we have \textcolor{black}{built} an approximation $X \approx H \Xi^T$ where $H$ is $n_t \times d$ score matrix such that $h^{(t)}$ is the $t^{th}$ column of $H$ and $\Xi$ the loading matrix. This is the  $d$-dimensional approximation of $X$ in  $\mathbb{R}^{n}$ that minimizes the mean squared error. Therefore, we have  $X G_*  \approx   H \Xi^T G(\Xi^T G)^{-1} = H$. 
Then $G_* \in \mathbb{R}^{n \times d}$ is the projection matrix from $X$ in the initial space to  $H$ in the reduced space and $\Xi^T$ is its reciprocal such that $ X G_* \Xi^T \approx H \Xi^T \approx X$.
It follows that, for a given reduced dimension $t \in \{ 1, \ldots, d \}$, a given point $x^r$ can be expressed in the original space along $t$:
\begin{equation} 
\begin{split}
& F_t : \mathbb{R}^{n } \xrightarrow[]{} \mathbb{R}^{n }  \\
& x^r  \mapsto \left[   \left[G_{*}\right]_{1}^{t} \left[x^r\right]_1 , \ldots, \left[G_{*}\right]_{n}^{t} \left[x^r\right]_{n}   \right].
\end{split}
\label{SMO_eq:KPLS_rota}
\end{equation}
\textcolor{black}{With PLS, we built a low-rank approximation and, in the following section, we show how to build the GP model in a small subspace instead of building it in the full space. 
The objective is twofold when the dimension increases. First, optimizing a small number of hyperparameters is much faster because, with more than 20 variables, building a GP is really prohibitive in terms of computational cost~\cite{Bouhlel18}. 
Second, optimizing the hyperparameters is harder and the resulting GP is often non-optimal and, not only is more costly than the KPLS model but also leads to some degraded performance~\cite{SciTech_cat}.
These two reasons motivated the need for the KPLS model described below.}
\subsubsection{KPLS for continuous variables}
\label{SMO_subsec:KPLS_cont}

The construction of the correlation matrix $R^{cont}(\theta^{cont})$ for continuous inputs, based on square exponential kernels (or Gaussian kernel~\cite{williams2006gaussian}) with PLS can be described as follows.  For a couple of continuous inputs $x^r\in \mathbb{R}^n$ and $x^s\in \mathbb{R}^n$, one sets:

% {\color{red}old version:
% \begin{equation}
% \begin{split}
%         [R^{cont}(\theta^{cont})]_{r,s} &= {\displaystyle \prod_{t=1}^{d} k_{t}^{cont} (F_t (x^r), F_t (x^s),\hat{\theta}_t^{cont}) } \\
%        &= {\displaystyle \prod_{t=1}^{d}{\displaystyle \prod_{j=1}^{n}  \exp \left(   - \hat{\theta}_t^{cont} \abs{g_{*j}^{(t)} x^r_j -g_{*j}^{(t)} x^s_j   }^p  \right)    }} \\
%         &= {\displaystyle \prod_{t=1}^{d}{\displaystyle \prod_{j=1}^{n}  \exp \left( - \hat{\theta}_t^{cont}   \left(g_{*j}^{(t)} \right)^p \abs{x^r_j - x^s_j}^p       \right) }} \\   
%        &=  \exp {\displaystyle \sum_{t=1}^{d}{\displaystyle \sum_{j=1}^{n}    - \hat{\theta}_t^{cont}   \left(g_{*j}^{(t)} \right)^p \abs{x^r_j - x^s_j}^p    }} \\  
%       &=  \exp {\displaystyle \sum_{j=1}^{n}    -   \lambda_j^{cont} \abs{x^r_j - x^s_j}^p      }      \\ 
%       &=   {\displaystyle \prod_{j=1}^{n}  \exp  \left( -   \lambda_j^{cont}  \abs{x^r_j - x^s_j}^p  \right)   }  ,    \\ 
% \end{split}        
% \label{SMO_eq:KPLS_exp_decomp}
% \end{equation}
% where $\lambda_j^{cont} = {\displaystyle \sum_{t=1}^{d}   \left( {g_{*j}^{(t)}} \right)^p \hat{\theta}_t^{cont}  }   $}

% {\color{blue}Ysf (new version): Paul check your implementation, an absolue value was missing in the expresion of $\lambda_j^{cont}$.
\begin{equation}
\begin{split}
        [R^{cont}(\theta^{cont})]_{r,s} &= {\displaystyle \prod_{t=1}^{d} k_{t}^{cont} (F_t (x^r), F_t (x^s),\hat{\theta}_t^{cont}) }  \\
       &= {\displaystyle \prod_{t=1}^{d}{\displaystyle \prod_{j=1}^{n}  \exp \left(   - \hat{\theta}_t^{cont} \left( \left[G_{*} \right]_j^{t} \left[x^r\right]_j -\left[G_{*} \right]_j^{t} \left[x^s\right]_j   \right)^2  \right)    }}\\ 
      &=   {\displaystyle \prod_{j=1}^{n}  \exp  \left( -   \theta_j^{cont}  \left(\left[x^r\right]_j - \left[x^s\right]_j \right)^2 \right)   }  ,    \\ 
\end{split}        
\label{SMO_eq:KPLS_exp_decomp}
\end{equation}
where $\theta_j^{cont} = {\displaystyle \sum_{t=1}^{d}   \left(  \left[ {G_{*}} \right]_j^t \right)^2 \hat{\theta}_t^{cont}  }$.
%}
%is the PLS approximation of $\theta_j^{cont}$.
%{\color{red}Ysf: All over the paper, you use $p=2$. why you decided here to have $p \in (0,2]$.} 
Clearly, in the continuous case, constructing $R^{cont}(\theta^{cont})$ would require the estimation of $d$ non-negative hyperparameters, $\textit{i.e.}$, $\theta^{cont} \in \mathbb{R}^d_{+}$, $d \ll n$.

%{\color{red} Ysf: Paul take you time to write Section 3.2 in a clean, clear and logical way. The section is the paper core. }
\subsection{Extension of PLS to matrix inputs with application to mixed-categorical GP}

This part presents the extension of PLS for a general categorical \textcolor{black}{GP} kernel. More precisely, in Section~\ref{SMO_sec:PLS_mat} we extend the PLS regression for matrix inputs and, in Section~\ref{SMO_subsec:KPLS_cat} we applied it to the \textcolor{black}{GP} kernels for categorical variables.

\subsubsection{A PLS framework for matrix inputs} % 
\label{SMO_sec:PLS_mat}

We consider a general categorical variable $c$ that can take $L$ different levels. In that context, we want to find a small $\ell \times \ell $ matrix $\hat{\Theta}$ to represent a bigger $L \times L$ correlation matrix $\Theta$, with $\ell < L$.  %\\
Recall that, from PLS, $G_*$ can be seen as the rotation matrix from the initial space to the reduced space~\cite{EGORSE}. By taking into account the symmetry of correlation matrices and their unit diagonal, we need to build a rotation matrix $G_*$ of dimension $\left(\frac{L(L-1)}{2} \times \frac{\ell(\ell-1)}{2} \right)$. 
The input dimension is denoted $D_{\inp}= \frac{L(L-1)}{2}$ and the output dimension is denoted $D_{\outp} = \frac{\ell(\ell-1)}{2}$. 

First, we want to construct the matrix $G_*$ aforementioned. %\\
For a given input $c^r$, its natural one-hot encoding  $e_{c^r}$ is a basis vector of dimension $L$~\cite{Mixed_Paul}. Meanwhile, the input that we need for $G_*$ is of dimension $D_{\inp}$ so, in order to have a vector data fitting the dimension, we propose to use a novel one-hot encoding relaxation that adds a new dimension for every cross-correlation term. Subsequently, the relaxation $\zeta_{c^r}$ is such that $\zeta_{c^r} \in \{0,1\}^{D_{\inp}}$: the $L-1$ terms that correspond to the correlation with the level taken by the input $c^r$ equal $1$ and all other terms take value $0$ and we call this relaxation "cross-levels encoding".
One can observe that the Hadamard product $\zeta_{c^r} \odot \zeta_{c^s} = 1 $  in the dimension corresponding to the correlation between the levels taken by ${c^r}$ and ${c^s}$ and zero everywhere else which is the property we were seeking for. In other words, for all $j \in {1, \ldots, D_{\inp}}, \quad$
\begin{equation}
\label{eq:prod_zeta_smo}
\begin{cases}
     [\zeta_{c^r} \odot \zeta_{c^s}]_j = 1, \quad \text{if } \zeta_{c^r_j} =1  \mbox{ and }  \zeta_{c^s_j} =1,  \\  
     [\zeta_{c^r} \odot \zeta_{c^s}]_j = 0, \quad \text{otherwise.}
\end{cases}
\end{equation}
Example~\ref{example:1} illustrates how the relaxed vectors $\zeta_{c^r}$ and $\zeta_{c^r}$ are built using a  simple use-case. 
\begin{example}\label{example:1}
Consider a categorical variable $c$ taking values in a color set of $L=4$ levels such that, for any point %{\color{red} Ysf: this  how we use cotes inside the math mode $``blue"$} 
$r$, $c^r \in \{ ``green", ``red", ``blue", ``yellow" \} $.
We want to represent the value of $c^r$ as $\zeta_{c^r} \in \{0,1\}^{D_{\inp}}$. 
In this case, ${D_{\inp}}$ = 6 which corresponds to the 6 possible correlations ( "blue-red", "blue-green", "blue-yellow", "red-green", "red-yellow" and "green-yellow"). 
For instance, if $n_t=3$ points are considered such that $(c^1, c^2, c^3) = \{ ``blue",``red", ``red" \}$, the first point $c^1 = ``blue"$ will be represented as $\zeta_{c^1} = (1,1,1,0,0,0)$, taking 1 for the dimensions related to the correlations involving "blue" and 0 everywhere else. Similarly, the second point $c^2 = ``red"$ will be represented as $\zeta_{c^2} = (1,0,0,1,1,0)$. And when taking the Hadamard product, $\zeta_{c^1} \odot \zeta_{c^2} = (1,0,0,0,0,0)$, the only dimension that takes value 1 corresponds exactly to the dimension representing "blue-red".
\end{example}

%\bigbreak
Using the relaxed DoE  $X = \{\zeta_{c^1} , \ldots, \zeta_{c^{n_t}} \}$ of dimension $ D_{\inp} \times n_t$, we can compute the rotation matrix $G_*$ of dimension  $ D_{\inp} \times D_{\outp} $ as in~\eqnref{SMO_eq:catKPLS}.
Our goal is to use the matrix $G_*$ to express an $L \times L$ matrix from a smaller $ \ell \times \ell$ matrix.  
%
%Denoting $ \left[G_{*} \right]_j^{j'}$ the $j,j'$ coefficient of $G_*$, 
The $L \times L $ symmetric matrix $\Theta$ with unit diagonal can be estimated using a smaller $ \ell \times \ell $ SPD matrix such that, for all $ j<j'$, one has  
%\vspace{-0.2cm}
\begin{equation}
\begin{split}
 \ & [\Theta]_{j}^{j'} =  {\displaystyle \sum_{t=1}^{\ell}{\displaystyle \sum_{t'=t+1}^{\ell}  \left( \left[ G_{*} \right]_{\psi(j,j',L)}^{\psi(t,t',\ell)} \right)^2 \ [\hat{\Theta}]_{t}^{t'}  }}, \\
%\forall j' > j,&  [\Lambda]_{j',j} = [\Lambda]_{j,j'}, \\
%\forall j=j',& [\Lambda]_{j,j'}=1.  \\
\end{split}    
\label{SMO_eq:true_pls_matrix_reduction}
\end{equation}
where we rely on a matrix-to-vector lexicographical transformation $\psi$ to insure that both the input vector of size $ \frac{L(L-1)}{2}$ and the output vector of size $\frac{ \ell (\ell-1)}{2}$ are valid representations of SPD matrices. For a given integer $n_\text{lev}$ and, for all $k \in \{1, \ldots, n_\text{lev} \}$ and $  k' \in \{k+1, \ldots, n_\text{lev} \}$, the mapping $\psi$ is given by:
%\vspace{-0.2cm}
\begin{equation}
\begin{split}
&\psi(k,k',n_\text{lev}) = \frac{1}{2} (  (n_\text{lev}-1)(n_\text{lev}-2) - (n_\text{lev}-k)(n_\text{lev}-k-1) )  +k'-1 . 
\end{split}
\label{SMO_eq:matrix_to_vector}   
\end{equation}
This formulation gives a sparse vector representation of the hyperparameters used to build the matrix  $\Theta$ by lexicographic order of the triangular superior part of the matrix. 
Notwithstanding, as we assumed that $\Theta$ is a symmetric matrix with unit diagonal, we could have defined $ [\Theta]^j_{j'}$ in~\eqnref{SMO_eq:true_pls_matrix_reduction} by the triangular inferior values. This would have led to the exact same kernel as what as been presented but with a slightly different definition of $\psi$.
With the expression of~\eqnref{SMO_eq:true_pls_matrix_reduction} we achieved to build a PLS approximation for matrices as intended. To finish with, we insure that~\eqnref{SMO_eq:true_pls_matrix_reduction} works for SPD matrices in order to build our GP upon it as described in the following section.

\begin{theorem}
\label{SMO_th:eq_gd_cr}
Assuming that all the entries of $\hat{\Theta}$ are in $[-1,1]$ and that $G^*$ is computed using PLS as in~\eqnref{SMO_eq:catKPLS}, the matrix $\Theta$ given by~\eqnref{SMO_eq:true_pls_matrix_reduction} also takes values in [-1,1].
\end{theorem}
\begin{proof}
Indeed, $G_*$ is a rotation matrix. Thus, for all $j\in\{1, \ldots,L\}$ and $  j' \in \{j+1, \ldots, L \}$, 
the $\psi(j,j',L)$-th row of $G_*$, given by $[G_*]_{\psi(j,j',L)} = \left\{ [G_*]_{\psi(j,j',L)}^{\psi(t,t',\ell)} \right\}_{\tiny \begin{array}{c}
1 \leq t \leq \ell, \\ \ t+1  \leq t' \leq \ell\end{array}}$, satisfies 
$$
\left([G_*]_{\psi(j,j',L)} \right)^\top \left([G_*]_{\psi(j,j',L)} \right) = \sum_{t=1}^{\ell}{\displaystyle \sum_{t'=t+1}^{\ell}  \left( \left[ G_{*} \right]_{\psi(j,j',L)}^{\psi(t,t',\ell)} \right)^2}=1.
$$
Hence, knowing that $\left|[\hat{\Theta}]_{t}^{t'}\right| \leq 1 $ for all $t,t'$, one has
\begin{eqnarray*}
   |[\Theta]_{j}^{j'} |&\le & {\displaystyle \sum_{t=1}^{\ell}{\displaystyle \sum_{t'=t+1}^{\ell}  \left( \left[ G_{*} \right]_{\psi(j,j',L)}^{\psi(t,t',\ell)} \right)^2 \ \left|[\hat{\Theta}]_{t}^{t'} \right| }}\le 1.
\end{eqnarray*}
\end{proof}
The matrix $\Theta$ is serving as a correlation matrix. For this purpose, it is essential that the matrix be SPD. To ensure its SPD nature, we check in our implementation if all of its eigenvalues are positive. If any eigenvalues are found to be negative, a nugget term is added to the covariance matrix to enforce the SPD property of the matrix $\Theta$. The nugget term allows us to mitigate numerical issues and maintain positive definiteness. It is worth noting that in all our numerical tests, the matrix $\Theta$ has been shown to be SPD. This suggests that if $\hat{\Theta}$ is SPD, then $\hat{\Theta}$ remains SPD, as discussed. Such a result seems not trivial to prove using ~\eqnref{SMO_eq:true_pls_matrix_reduction}.
%
%
%If the correlation matrix $\Theta$ is not SPD, it is possible to add a nugget to the correlation matrix $R$ between the inputs. Notwithstanding, we used no nugget in the following which indicates that the defined matrix may be SPD. 
%
Note also that $\Theta$ gives a good approximation of the correlation matrix between the levels of a categorical variable.
This approximation can be used to understand the structure of our modeling problem as in Section~\ref{SMO_sec:Results}.

% \begin{theorem}
% \label{SMO_th:eq_gd_cr}
% $[\Theta]$ constructed as in~\eqnref{SMO_eq:true_pls_matrix_reduction} is symmetric positive definite with values in $[-1,1]$. 
% \end{theorem}

% \begin{proof}
% By construction, we have $\forall j, j',  [\Theta]_{j,j'} = [\Theta]_{j',j}$ so $[\Theta]$ is symmetric.
% The construction implies the so-called rotation vectors $g_*$ as defined in Section~\ref{SMO_subsec:PLS_cont}. Orthogonal projections are positive definite operators and change of basis preserves eigenvalues. We are using the PLS regression over vectors thanks to our lexicographic order $\Psi$ and we know from the continuous case that the PLS regression defines a positive definite operation over vectors~\cite{PLS_Li_2002,Bouhlel18}. Also, the PLS does not change the values the hyperparameters can take as it is an unitary transform. As $[\hat{\Theta}_i]$ is constructed by hypersphere decomposition, the SPD~\cite{hypersphere} property is ensured with values in $[-1,1]$, and thus it is the case for $[\Theta_i]$.  
% \end{proof}

% \begin{corollary}

% Let define the matrix $ [\hat{\Theta}_i]$ of dimension $ L_i \times L_i $ with the homoscedastic hypersphere model of Zhou et al.~\cite{Pelamatti}. The higher dimension matrix $[\Theta_i]$  constructed as in~\eqnref{SMO_eq:true_pls_matrix_reduction} being SPD, we have also defined a reduced order model for every correlation matrix.
% \end{corollary}

\subsubsection{A new KPLS model for categorical variables}  
\label{SMO_subsec:KPLS_cat}

For a given categorical variable $i$, we want to express the matrix $R_i(\Theta_i)$ with less than  $D_{\inp}=\frac{L_i(L_i-1)}{2}$ hyperparameters $\Theta_i$. 
To do so, we generalize the KPLS method of Bouhlel~\textit{et al.}~\cite{Bouhlel18} for any correlation matrix. 
Let $\hat{\Theta}_i$ be a $\ell_i \times \ell_i$ SPD matrix defined on the reduced space whose values are in $[-1,1]$ constructed by homoscedastic hypersphere decomposition~\cite{Zhou}. The $D_{\outp} = \frac{\ell_i(\ell_i-1)}{2} $ correlation parameters of $\hat{\Theta}_i$ can be optimized by MLE from the scores projected data $H_i = X_i G_*$ as in the continuous case.  
%
%
% Recall that the new relaxation $\zeta$ is build such that $\zeta_{c_i^r} \odot \zeta_{c_i^s} = 1 $  in the dimension corresponding to the correlation between the levels taken by ${c^r}$ and ${c^s}$ and zero everywhere else. Based on this relaxation,  we can express the homogeneous kernel as follows. Namely, if $c^r_i =c^s_i$, $ [R_i(\Theta_i)]_{\ell^i_r,\ell^i_s} = 1$ and, if $ c^r_i \neq c^{s}_i$, we get
% \begin{equation}
% \begin{split}
%  [R_i(\Theta_i)]_{\ell^i_r,\ell^i_s} 
%  & = {\displaystyle \prod_{j=1}^{L_i} {\displaystyle \prod_{j' =j+1}^{L_i}
%         \kappa ( [ \Phi(\Theta_i) ]_{j,j'} \      ( \zeta_{c^r_i})_{\psi(j,j',L_i)}  (\zeta_{c^s_i})_{\psi(j,j',L_i)}    }}   ) \\
%         &= 
% \  \kappa ( [ \Phi(\Theta_i) ]_{{ \ell_i^r},{\ell_i^s}} ),  \\
% \end{split}    
% \label{SMO_eq:KRG_in_PLS_space}
% \end{equation}
% with $\psi$ defined in~\eqnref{SMO_eq:matrix_to_vector}. This expression is the same as in~\eqnref{SMO_eq:homo_HS} meaning that~\eqnref{SMO_eq:KRG_in_PLS_space} is the homogeneous kernel in our new formalism. 
%
%
Based on~\eqnref{SMO_eq:true_pls_matrix_reduction} for matrix PLS, 
%and~\eqnref{SMO_eq:KRG_in_PLS_space} for KPLS
we can introduce our new HH and EHH KPLS kernels depending only on $D_{\outp}$ hyperparameters defined as follows. 
Recall that the matrix $C(\Theta_i)\in \mathbb{R}^{L_i \times L_i}$ is lower triangular and built using a hypersphere decomposition~\cite{HS,HS_Jacobi} and that the variable $\epsilon$ is a small positive constant. 
\begin{itemize}
    \item  
\label{def:HH_PLS}
The HH KPLS kernel is given by $\kappa = \mathbb{I}_{L_i}$, $[\Phi(\Theta_i)]_{j}^{j} =1$ and, for all $j \neq j'$,
%\vspace{-0.25cm}
$$ [\Phi(\Theta_i)]^{j'}_{j}  = \frac{1}{2} \  {\displaystyle \sum_{t=1}^{\ell}{\displaystyle \sum_{t'=t+1}^{\ell} \left( \left[ G_{*} \right]_{\psi(j,j',L)}^{\psi(t,t',\ell)} \right)^2 \ [C(\Theta_i) C(\Theta_i)^\top]_{t}^{t'}  }}.$$
\item
\label{def:EHH_PLS}
The EHH KPLS kernel is given by $ \kappa(\phi) = \exp (-\phi) $, $[\Phi(\Theta_i)]_{j}^{j} =0$ and, for all $j \neq j'$,
%\vspace{-0.25cm}
$$ [\Phi(\Theta_i)]^{j'}_{j}  =  {\displaystyle \sum_{t=1}^{\ell}{\displaystyle \sum_{t'=t+1}^{\ell}   \left( \left[ G_{*} \right]_{\psi(j,j',L)}^{\psi(t,t',\ell)} \right)^2 \ \frac{\log \epsilon}{2} ([C(\Theta_i) C(\Theta_i)^\top]_{t}^{t'}-1) }}.$$
\end{itemize}
%$[\Theta_i] $ can be used as a starting point for a maximization of the likelihood function for a standard Kriging model on the whole starting space, this method called KPLS(+K) can be generalized for categorical kernels. 
%Note that the same categorical input can lead to different output ${y}$ depending on the continuous and discrete variables. Therefore, we need $\frac{\ell_i (\ell_i+1)}{2}+1$ hyperparameters, by taking into account the noise evaluation.
%
For EHH KPLS kernels, the proposed KPLS model as given by~\eqnref{SMO_eq:true_pls_matrix_reduction} can be shown as a natural extension of the continuous KPLS as proposed by Bouhlel~\textit{et al.}~\cite{Bouhlel18} (described also in Section~\ref{SMO_subsec:KPLS_cont}). The result is shown hereinafter.

%{\color{red} Ysf: Pas claire, Tu veux dire une extension pour KPLS dans le cas dles variables continues. Tu as besoin aussi de préciser que pour un noyaux exp. LEs hypothèses sur $\hat{\Theta}$ en entrée manquent aussi }
\begin{theorem}
For a correlation matrix $\hat{\Theta}_i$, the projection formula used in~\eqnref{SMO_eq:true_pls_matrix_reduction} extends the continuous KPLS to categorical matrices using an exponential kernel. 
\end{theorem}
%
% {\color{red}Ysf: write a clean latex. Avoid to mix text and equation. Paul, this is your third paper, you should know now how we write properly a paper. Look at the first paper, e.g., Eq (17).}
\begin{proof}
The KPLS kernel, for exponential kernel, is based on the fact that, for a given reduced dimension $t \in \{ 1, \ldots, d \}$, a given point $x^r$ can be expressed in the original space along $t$ as in~\eqnref{SMO_eq:KPLS_rota}. 
%
% \begin{equation} 
% \begin{split}
% & F_t : \mathbb{R}^{n } \xrightarrow[]{} \mathbb{R}^{n }  \\
% & x^r  \mapsto \left[   g_{*1}^{(t)} x^r_1 , \ldots, g_{*n}^{(t)} x^r_{n}   \right].
% \end{split}
% \end{equation}
%
Therefore, we apply the same transformation to our relaxed inputs and then apply the transformation $\psi$ to have a matrix formulation from the relaxed vectors. 
This leads to the natural way to express our new EHH KPLS kernel defined as $[R_i(\Theta_i)]_{\ell^r_i}^{\ell^s_i}   = 1$ if $ c^r_i= c^s_i$ and otherwise, %$[R_i(\Theta_i)]_{\ell^r_i}^{\ell^s_i} =$
\begin{equation*}
\begin{split} 
& [R_i(\Theta_i)]_{\ell^r_i}^{\ell^s_i}\\
 &= {\displaystyle \prod_{t=1}^{\ell_i} {\displaystyle \prod_{t' =t+1}^{\ell_i} {\displaystyle \prod_{j=1}^{L_i} {\displaystyle \prod_{j' =j+1}^{L_i}
        \exp \left[ - \left(   \left[G_* \right]_{\psi(j,j',L_i)}^{\psi(t,t',\ell_i)} [  \zeta_{c^r_i}]_{\psi(j,j',L_i)}  \ \left[G_* \right]_{\psi(j,j',L_i)}^{\psi(t,t',\ell_i)} [\zeta_{c^s_i}]_{\psi(j,j',L_i)} \right) \ [\hat{\Theta}_i]_{t}^{t'} \  \right] }}}}  \\
     &     = {\displaystyle \prod_{t=1}^{\ell_i} {\displaystyle \prod_{t' =t+1}^{\ell_i} {\displaystyle \prod_{j=1}^{L_i} {\displaystyle \prod_{j' =j+1}^{L_i}
        \exp \left[ - \left(    [  \zeta_{c^r_i}]_{\psi(j,j',L_i)}   [\zeta_{c^s_i}]_{\psi(j,j',L_i)} \right) \ \left( \left[G_* \right]_{\psi(j,j',L_i)}^{\psi(t,t',\ell_i)} \right)^2 [\hat{\Theta}_i]_{t}^{t'} \  \right] }}}}  \\
    &       =\exp \left[  {\displaystyle \sum_{j=1}^{L_i} {\displaystyle \sum_{j' =j+1}^{L_i} {\displaystyle \sum_{t=1}^{\ell_i} {\displaystyle \sum_{t' =t+1}^{\ell_i} - \left( \left[G_* \right]_{\psi(j,j',L_i)}^{\psi(t,t',\ell_i)} \right)^2  [\hat{\Theta}_i]_{t}^{t'} \   }} \ 
         \left(    [  \zeta_{c^r_i}]_{\psi(j,j',L_i)}   [\zeta_{c^s_i}]_{\psi(j,j',L_i)} \right) }} \right] \\
    &       = \exp \left[ {\displaystyle \sum_{j=1}^{L_i} {\displaystyle \sum_{j' =j+1}^{L_i}
        -  [{\Theta}_i]_{j}^{j'}  \left(  [\zeta_{c^r_i}]_{\psi(j,j',L_i)}   [\zeta_{c^s_i}]_{\psi(j,j',L_i)} \right) \ \   }} \right]  \\
       &       = \exp \left[ {\displaystyle \sum_{j=1}^{L_i} {\displaystyle \sum_{j' =j+1}^{L_i}
        -  [{\Theta}_i]_{j}^{j'}  \left(  \delta_{j,\ell^r_i} \delta_{j',\ell^s_i}  \right) \ \   }} \right]  \\
 &          = \exp \left[ -[{\Theta}_i]_{\ell_i^r}^{\ell_i^s}  \right], %\qed 
\end{split}    
\label{SMO_eq:KPLS_in_PLS_space}
\end{equation*}
where $\delta_{i,j}$ is the Kronecker symbol (i.e., $\delta_{i,i}=1$ and $\delta_{i,j}=0$ for all $i\neq j$) and $\Theta \in \mathbb{R}^{L_i \times L_i }$ is given by~\eqnref{SMO_eq:true_pls_matrix_reduction}.

% as
% \begin{equation}
% \begin{split}
%  \forall j< j', \ & [\Theta]_{j}^{j'} =  {\displaystyle \sum_{t=1}^{\ell}{\displaystyle \sum_{t'=t+1}^{\ell}  \left( \left[ G_{*} \right]_{\psi(j,j',L)}^{\psi(t,t',\ell)} \right)^2 \ [\hat{\Theta}]_{t}^{t'}  }}. \\
% %\forall j' > j,&  [\Lambda]_{j',j} = [\Lambda]_{j,j'}, \\
% %\forall j=j',& [\Lambda]_{j,j'}=1.  \\
% \end{split}    
% \end{equation}
% Moreover, the SPD matrix $\hat{\Theta}$ can defined using a smaller EHH kernel~\cite{Mixed_Paul} and we obtain the EHH KPLS kernel given by $\kappa(\phi) = \exp (-\phi)$, $[\Phi(\Theta_i)]_{j}^{j} =0$ and 
% %{\color{red} Ysf: LEs 2 sommes sont à l intérieurs du exp, non ?} 
% $$ [\Phi(\Theta_i)]_{j \neq j'}  = {\displaystyle \sum_{t=1}^{\ell}{\displaystyle \sum_{t'=t+1}^{\ell}   \left( \left[ G_{*} \right]_{\psi(j,j',L)}^{\psi(t,t',\ell)} \right)^2 \ \frac{\log \epsilon}{2} ([C(\Theta_i) C(\Theta_i)^\top]_{t}^{t'}-1) }}.$$
% %
% %
\end{proof}
% \begin{corollary}
% Let define the matrix $ \hat{\Theta}_i$ of dimension $ \ell_i \times \ell_i $ with the homoscedastic hypersphere model of Zhou et al.~\cite{Zhou}. The higher dimension matrix $\Theta_i$ constructed as in~\eqnref{SMO_eq:true_pls_matrix_reduction} being SPD, we have also defined a reduced order model for every correlation matrix based on EHH kernel~\cite{Mixed_Paul}.
% \end{corollary}
In the next section, we will see how to apply our new KPLS matrix-based GP on analytical and engineering problems. We will demonstrate how our surrogate models can provide insights into the underlying structure of the correlation matrix and how it can be utilized for \textcolor{black}{BO} when dealing with structural and multidisciplinary problems. 

\section{Results and discussion}
\label{SMO_sec:Results}

In this section, we demonstrate how our GP performs over various test cases and compare it to other GP models.
To begin with, Section~\ref{SMO_sec:imp_detail} gives the details of the implementation used for the following computer experiments.
Next, Section~\ref{SMO_subsec:model_val} illustrates the GP models on an analytical and on a structural test case.
Finally, section~\ref{SMO_sec:BO_val} presents the use of these GP models for \textcolor{black}{BO} with mixed variables on an analytical test case and then in the context of MDO for aircraft design. \textcolor{black}{The considered test cases in this section, as well as the number of hyperparameters related to each kernel, are listed in \tabref{SMO_tab:test_case}.}

\begin{table}[!h]
\centering
 \caption{\textcolor{black}{Number of variables and hyperparameters for the test cases in this study.}}
 \label{SMO_tab:test_case}
\small
\begin{tabular*}{\linewidth}{ccc}
\hline
\textbf{\textcolor{black}{ $ \ \ \quad \quad \quad $ Problem $\quad \quad \quad \ \ $ }} & \textcolor{black}{ \# of variables} & 
\begin{tabular}{ccc}
  \ & \textcolor{black}{ \# of hyperparameters } & \  \\
  \hline
   \textcolor{black}{ GD} & \textcolor{black}{ CR} &  \textcolor{black}{HH}  \\
\end{tabular} \\
\hline
\\
\end{tabular*}
\small
\vskip-10pt
\begin{tabular*}{\linewidth}{ccccc} 
\textbf{\textcolor{black}{Categorical cosine problem}} & $ \quad \quad \quad $ \textcolor{black}{2} & $ \quad \quad \ \    $ \textcolor{black}{2} & $ \quad \quad \quad \quad \  $  \textcolor{black}{14} & $ \quad  \quad  \quad \quad \  $  \textcolor{black}{79}  \\   
\textbf{\textcolor{black}{Cantilever beam}} & $ \quad \quad \quad $ \textcolor{black}{3} &  $ \quad \quad \ \   $ \textcolor{black}{3} & $ \quad \quad \quad \quad \ $ \textcolor{black}{14} & $ \quad  \quad  \quad \quad \  $ \textcolor{black}{68}  \\   
\textbf{\textcolor{black}{Toy function}} & $ \quad \quad \quad $ \textcolor{black}{2} & $ \quad \quad \ \   $ \textcolor{black}{2} & $ \quad \quad \quad \quad \ $ \textcolor{black}{11} & $ \quad  \quad  \quad \quad \  $ \textcolor{black}{46}  \\   
\textbf{\textcolor{black}{\texttt{DRAGON} aircraft concept}} & $ \quad \quad \quad $ \textcolor{black}{12} & $ \quad \quad \ \   $ \textcolor{black}{12} & $ \quad \quad \quad  \quad  \ $ \textcolor{black}{29} &  $ \quad  \quad  \quad \quad \  $ {\textcolor{black}{147}}  \\   
\hline
\end{tabular*}
\end{table}

\subsection{Implementation details}
\label{SMO_sec:imp_detail}
Optimizing the likelihood with respect to the hyperparameters necessitates the use of an efficient gradient-free algorithm. In this study, we have employed COBYLA~\cite{COBYLA} from the Python library Scipy, which employs default termination criteria related to the trust region size. Since COBYLA is a local search algorithm, we have employed a multi-start technique for enhanced robustness. Our models and their implementation can be accessed in the SMT v2.0 toolbox~\cite{SMT2019, saves2023smt}. In SMT 2.0\footnote{\url{https://smt.readthedocs.io/en/latest/}}, the default number of starting points for COBYLA is set to 10, distributed evenly. We utilize a straightforward noiseless Kriging model with a constant prior. It is important to note that the absolute exponential kernel and the squared exponential kernel behave similarly for categorical variables. The correlation values range between $2.06e-9$ and $0.999999$ for both continuous and categorical hyperparameters. Consequently, we select the constant $\epsilon$ to correspond to a correlation value of $2.06e-9$. The DoE are generated through Latin Hypercube Sampling (LHS)~\cite{LHS}, and the validation sets consist of evenly spaced points.

For \textcolor{black}{BO} without constraint, we are using the EGO method of SMT 2.0 with the aforementioned GP models from the same toolbox. For \textcolor{black}{BO} under constraints, we are using the SEGOMOE method~\cite{bartoli:hal-02149236} with the mean criterion for the  metamodels of constraints~\cite{SEGO-UTB}. We are using the same GP models for both objective and constraints. The optimization of the WB2s infill criterion~\cite{bartoli:hal-02149236} is done using SNOPT~\cite{gill2005snopt}.
To compare \textcolor{black}{BO} with, we used a Multi-Objective Evolutionary Algorithm (MOEA)~\cite{Petrowski2017} named the Non-dominated Sorting Genetic Algorithm II (NSGA-II) \cite{nsga2} due to its low configuration effort and high performance. The NSGA-II algorithm used is the implementation from the toolbox pymoo~\cite{pymoo} with the default parameters (probability of crossover = 1, eta = 3). Pareto fronts are not relevant in our study as we are considering single-objective optimization. \textcolor{black}{We note that, although NSGA-II is designed for multi-objective optimization, for the purpose of establishing a baseline reference (for comparison), we have used NSGA-II to solve our mono-objective optimization problem. In fact, to the best of our knowledge, NSGA-II is the only open-source optimization solver available for addressing mixed-variable constrained optimization.}
\textcolor{black}{In this paper, all results are obtained using an Intel® Xeon® CPU E5-2650 v4 @ 2.20 GHz core and 128 GB of memory with a Broadwell-generation processor front-end and a compute node of a peak power of 844 GFlops.}

\textcolor{black}{Note that when using KPLS, the GP models built based on the EHH and HH kernels demonstrate comparable performance, with a slight advantage for HH in our numerical tests. For this reason, and to enhance the readability of the numerical section, we have decided to report only results on the HH kernel. Nevertheless, the method related to the use of KPLS within the EHH kernel (i.e., EHH-PLS) is available in the SMT 2.0 toolbox.}

\subsection{Surrogate modeling}
\label{SMO_subsec:model_val}

In this section, we validate our model on both a state-of-the-art modeling problem in Section~\ref{SMO_subsec:catcos} and on an structural cantilever beam problem in Section~\ref{SMO_subsec:beam}. More precisely, the matrix based PLS model is compared with literature models like  GD, CR or HH. This section shows that the PLS information can capture the shape of the correlation matrix between the various levels of a categorical variable.

\subsubsection{Analytic validation on a categorical cosine problem \textcolor{black}{($n= 1$, $m=0$, $l=1$ and $L_1=13$)}}
\label{SMO_subsec:catcos}

In this section, we investigate the categorical cosine problem, as outlined in~\cite{Roustant}, to showcase the behavior of the proposed kernels. The black-box function, denoted as $f$, relies on both a continuous variable within the range of $[0,1]$ and on a categorical variable with 13 distinct levels. 
\textcolor{black}{Consequently, the relaxed dimension(i.e., the number of hyperparameters) is 14 for the construction of a continuous GP with CR and the most general GP with our new relaxation is of dimension 79.} 

Appendix~\ref{SMO_subsec:cosine} provides a detailed description of this function.
A given point of the DoE is denoted as $w= (x,c)$, where $x$ represents the continuous variable and $c$ represents the categorical variable. This work aims at modeling the interplay between the various variables (together with their respective levels) as well as their impact on the targeted function. Notably, the categorical variable exhibits two distinct groups of curves, each comprising a subset of the 13 levels. The first group encompasses levels 1 to 9, while the second group consists of levels 10 to 13. Within each group, we observe strong positive correlations, implying that variables within the same group exhibit a similar behavior. Conversely, strong negative correlations manifest between the two groups, indicating distinct behavior and characteristics between them.  
In this example, the number of relaxed dimensions for continuous relaxation is 14. A LHS DoE with 98 points ($14\times 7$, if 7 points per dimension are considered) is chosen to built the GP models. On this test case, the number of hyperparameters to optimize is therefore $2$ for GD and HH with 2D PLS, $14$ for CR and $79$ for HH as indicated in~\tabref{SMO_tab:resRoustant}. 

\begin{figure}[!h]
\begin{center}
\vspace{-0.25cm}
        \subfloat[GD kernel (2 hyperparameters: 1 cat. and 1 cont.)]{
      \centering
		\includegraphics[clip=true, height=4.5cm, width=5cm]{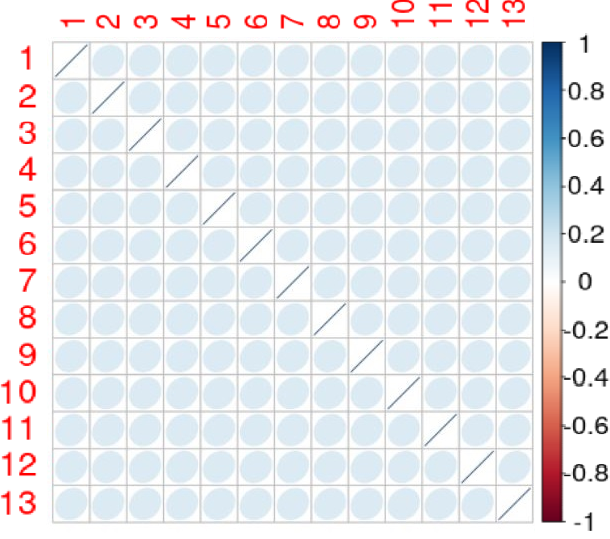}
    	\includegraphics[clip=true, height=4.5cm, width=8.6cm]{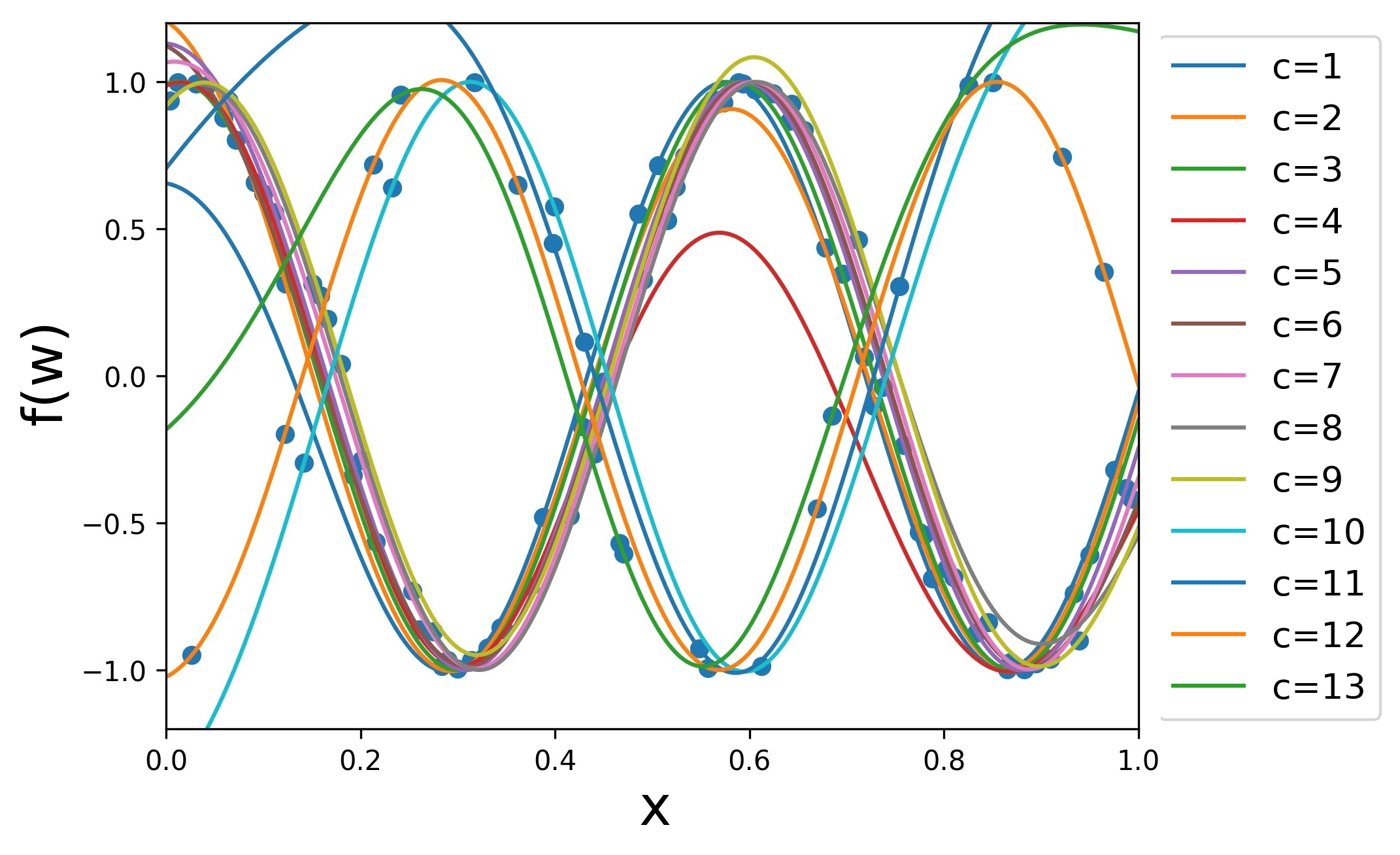}
    }
\vspace{-0.25cm}
      \subfloat[HH with PLS kernel (2 hyperparameters: 1 cat. and 1 cont.)]{
      \centering
		\includegraphics[clip=true, height=4.5cm, width=5cm]{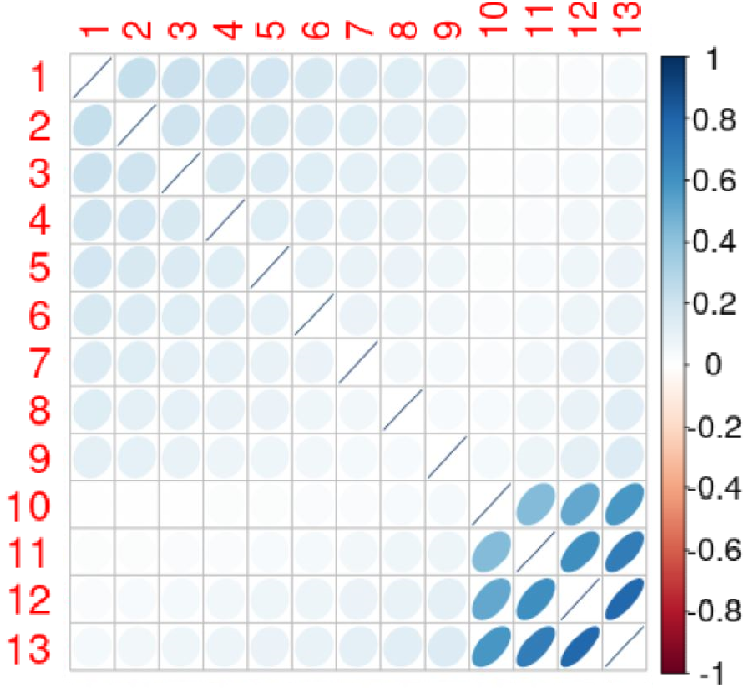}
    	\includegraphics[clip=true, height=4.5cm, width=8.6cm]{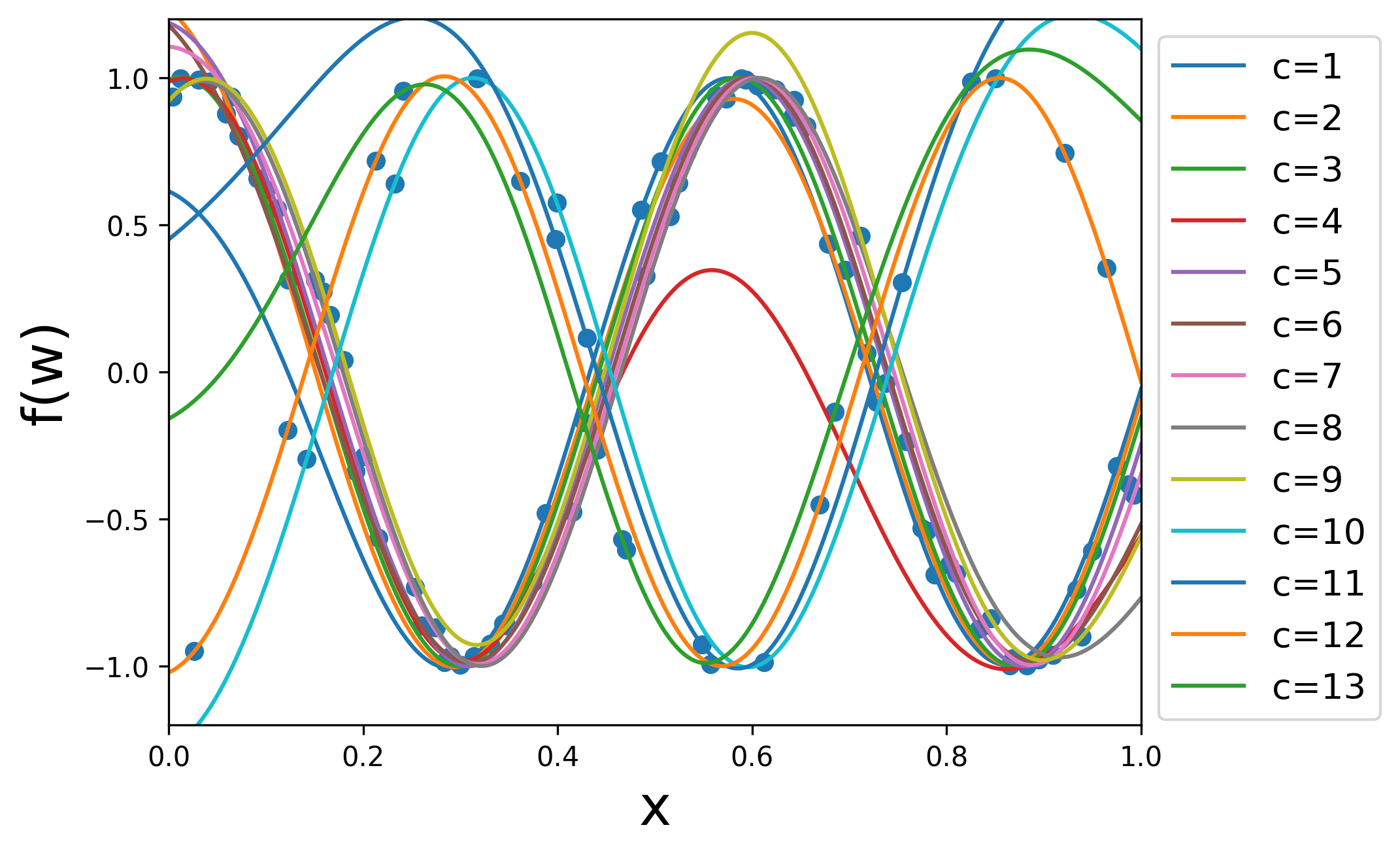}
    }
\vspace{-0.25cm}
      \subfloat[CR kernel (14 hyperparameters: 13 cat. and 1 cont.)]{
      \centering
		\includegraphics[clip=true, height=4.5cm, width=5cm]{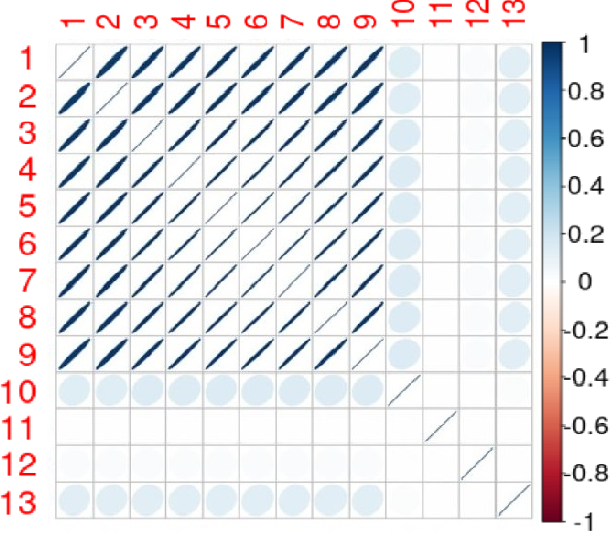}
    	\includegraphics[clip=true, height=4.5cm, width=8.6cm]{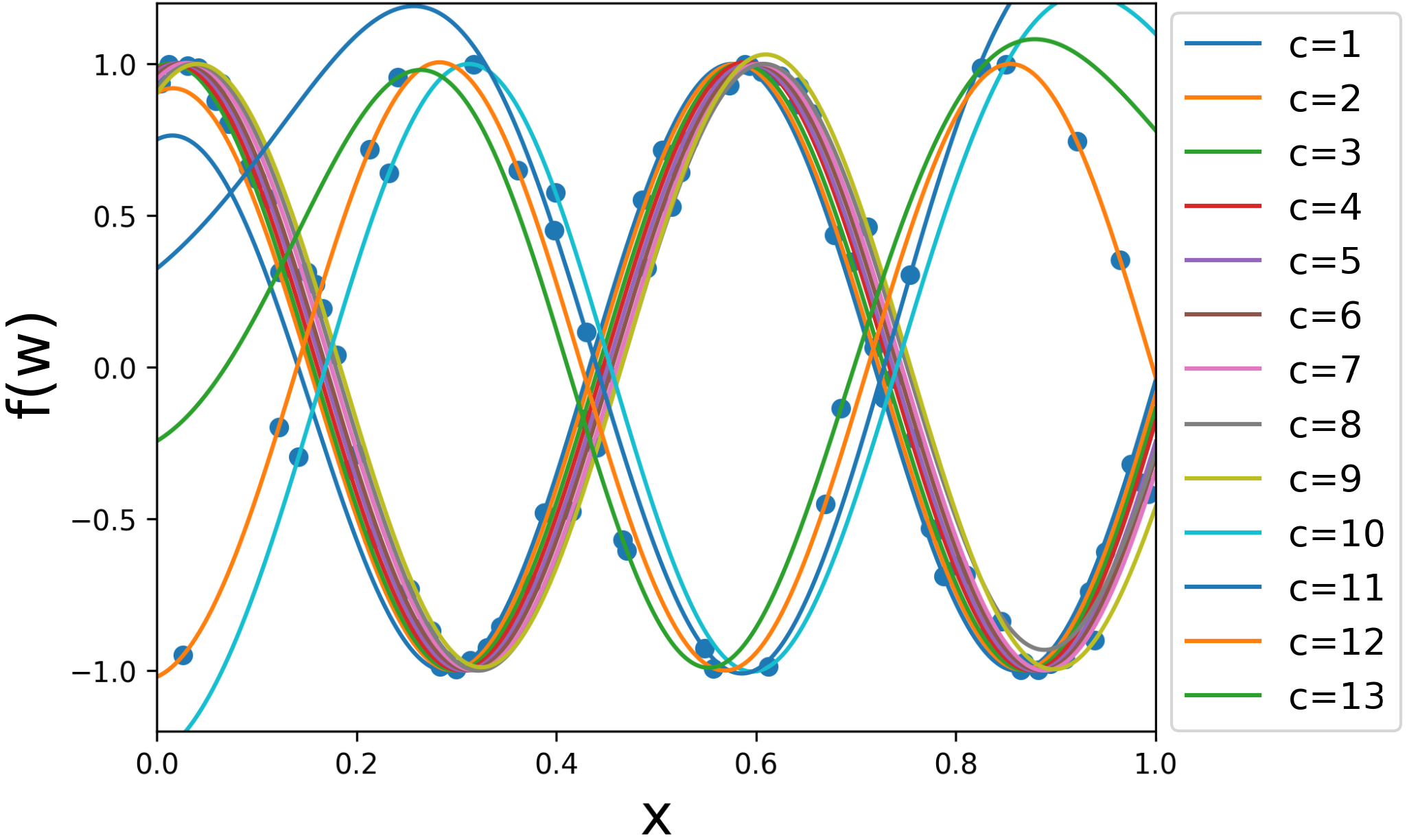}
    }
\vspace{-0.25cm}
  %       \subfloat[EHH kernel (79 hyperparameters) ]{
  %     \centering
		% \includegraphics[clip=true, height=4.5cm, width=5cm]{images/corr_homo_R.png}
  %   	\includegraphics[clip=true, height=4.5cm, width=8.6cm]{images/HOMO_50_Roustant_curves.png}
  %   \label{SMO_Roustant_hs_gsh}
  %   }
    \subfloat[HH kernel (79 hyperparameters: 78 cat. and 1 cont.)]{
  \centering
	\includegraphics[clip=true, height=4.5cm, width=5cm]{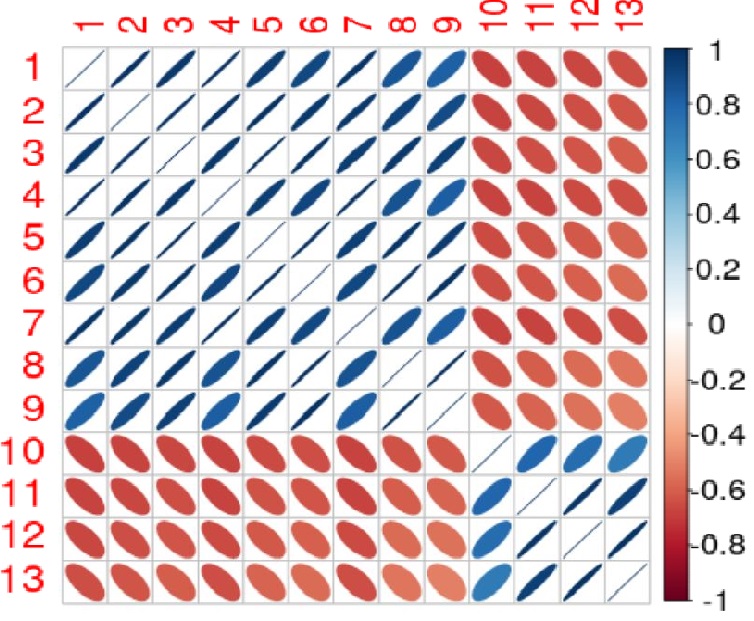}
	\includegraphics[clip=true, height=4.5cm, width=8.6cm]{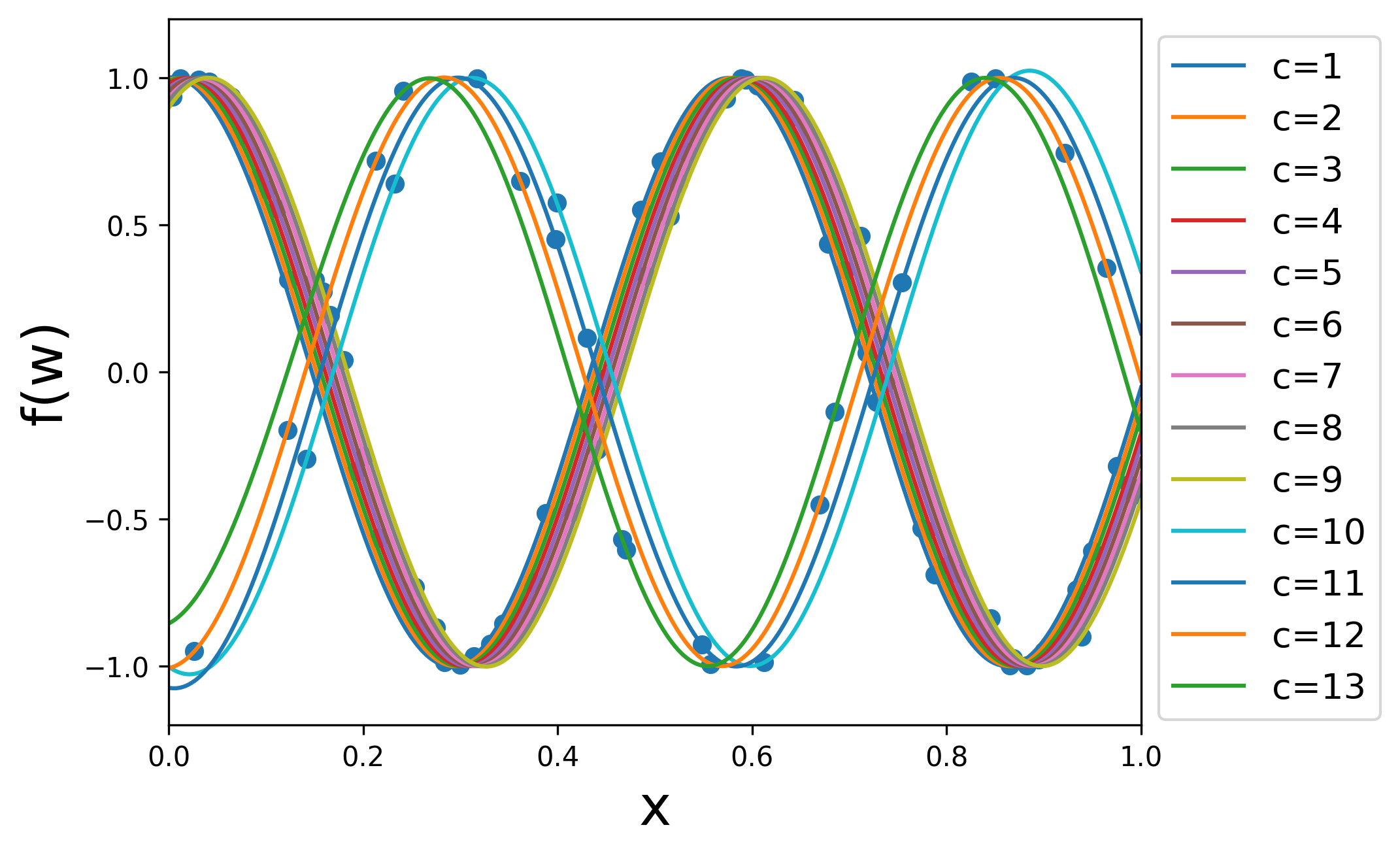} 
     }
\caption{Correlation matrices and associated predictions on the cosine problem using a DoE of 98 points. }
\label{SMO_Roustant_comp}
\end{center} 
\end{figure}

For GD, CR, HH and HH with PLS, the associated mean posterior models are shown on~\figref{SMO_Roustant_comp} on the right and the estimated correlation matrices $R_i=R_1^{cat}$ are displayed on the left.
The latter matrices can be interpreted as such: 
for two given levels $\{\ell_r^1,\ell_s^1\}$, the correlation term $[R_1]_{\ell_r^1,\ell_s^1}$ is in blue for correlation value close to 1, in white for correlations close to 0 and in red for value close to -1; moreover the thinner the ellipse, the higher the correlation and we can see that the correlation between a level and itself is always 1. 
At first glance, one can see, on~\figref{SMO_Roustant_comp}, that the predicted values remain properly within the interval $[-1,1]$ only with the HH kernel and looked badly estimated with GD. To quantify this assumption, we compute the Root Mean Square Error (RMSE) and Predictive Variance Adequacy (PVA)~\cite{PVA} for every model as 
\begin{equation*}
\label{SMO_eq:RMSE}
\mbox{RMSE}= \sqrt{\underset{i=1}{\overset{N}{\sum}} \frac{1}{N} \left( \hat{f}(w_i) - f(w_i) \right)^2} \quad\mbox{and} \quad 
\mbox{PVA}= \log \left( {\underset{i=1}{\overset{N}{\sum}} \frac{1}{N} \frac{\left( \hat{f}(w_i) - f(w_i) \right)^2}{ [\sigma^f(w_i)]^2 }  } \right),
\end{equation*}
where $N$ is the size of the validation set, $\hat{f}(w_i)$ and $\sigma^f(w_i) $ are the mean and standard deviation predictions of our GP model at a point $w_i$, $f(w_i)$ is the associated true value and the validation set consists of $N=13000$ evenly spaced points (see Appendix~\ref{SMO_subsec:cosine}). The values are reported in~\tabref{SMO_tab:resRoustant} and show that the PVA is constant indicating that the estimation of the variance is proportional to the RMSE. Nevertheless, the PVA is smaller for the EHH and HH kernels because the optimization has not converged yet after 887 seconds (79 parameters being hard to optimize in that case). However, a longer run gives a RMSE of 1.280 and a PVA of 21.95 for HH, which, once again is around the other PVA values. Note that the RMSE obtained with HH and EHH are significantly smaller than the errors obtained with all other methods, even with an incomplete optimization.

In~\tabref{SMO_tab:resRoustant}, we also included the results obtained for different sizes of the PLS reduced matrix, namely $2 \times 2,3 \times 3,4 \times 4$ and $5 \times 5$. As expected, these results show that the more complex the model, the smaller the RMSE. Moreover, with the same number of hyperparameters the HH 2D PLS outperforms GD in terms of accuracy and retrieves part of the correlation matrix structure. Also, the HH 2D PLS model can lead to almost similar performance as CR with significantly less hyperparameters.
\textcolor{black}{Indeed, PLS reduces slightly the accuracy of the model because KPLS is a simplified model optimized in a small subspace but it reduces the run time by a factor of around 10 on this categorical cosine test case. Both the CR and HH models have been built for illustration purposes but are intractable in time for real test cases and the latter is the reason for developing reduced order model in this paper.}
Notwithstanding, in terms of computational costs, computing PLS for matrices could be time-consuming (130 seconds for 2D PLS versus 2.5 seconds for GD). However, this cost is a fixed cost making our method of particular interest for large datasets and high number of dimensions.
%whereas the computational time scales poorly with the number of hyperparameters and the size of the dataset
%
%\vspace{-0.25cm}
\begin{table}[!h]
\centering
 \caption{Kernel comparison for the cosine test case in terms of number of hyperparameters, time, RMSE and PVA metrics.}
 \label{SMO_tab:resRoustant}
\small
\begin{tabular}{cccccc}
\hline
  \textbf{Kernel} &  \# of Hyperparam. & run time (s) & RMSE & PVA   \\
  \hline 
  \textbf{GD} &  2 & 2.5 & 30.079&  21.99  \\   
  %HH+PLS & 4 & 27.284  & 21.95  \\
 \hdashline 
  \textbf{CR} &  14 & 19& 22.347 & 23.04 \\
  \textbf{CR+PLS} &   2 & 2.8 & 26.376  & 21.87  \\
  \hdashline 
%  HH  & 79 & 1.280   & 24.31  \\
  \textbf{HH}  & 79 & 887 & 5.330   & 15.34  \\
 \textbf{HH+PLS(2x2)} & 2 & 130 & 26.087  & 21.86  \\
  \textbf{HH+PLS(3x3)} & 4 & 326 & 25.504  & 21.95  \\
  \textbf{HH+PLS(4x4)} & 7 & 819& 23.01 & 22.13  \\
  \textbf{HH+PLS(5x5)} & 11 & 1787& 23.04 & 22.13  \\
  \hdashline 
  \textbf{EHH}  & 79 & 959 & 6.858   & 15.46  \\
 %   FE & 92 & 22.610 & 642.0 \\
\hline
\end{tabular}
\end{table}
%
%
%\vspace{-0.25cm}
 We also include the metrics for CR+PLS but this method combines all hyperparameters and does a unique PLS in a continuous space, therefore, this method cannot be used to retrieve the correlation matrix and to study the categorical variable, it has no explainability. Nevertheless, in terms of predictive power, this method is associated with a computational time similar to GD for a better predictive power which makes it a good method for quick prediction on uncertain zones as in the context of \textcolor{black}{BO}. For the method developed in this paper, we tested several PLS components but these components correspond to a small correlation matrix and therefore should be of the form $\frac{n(n+1)}{2}, n \in \mathbb{N}$. We computed the models with the first number of PLS components ($1,3,6,10$) as indicated in~\tabref{SMO_tab:resRoustant} but the predictive gain is not significant, especially given the increase in run time to build the model when the number of components increases. In particular, after 6 components (4x4 matrix), we have a cost comparable with the full model, and, even if these 6 hyperparameters are optimized completely, the prediction is still rough whereas the incomplete optimization of the HH model gives a significantly better prediction still. 
For that reason, in the following experiments and comparisons, we will stick to the 1 component (2x2 matrix) PLS to retrieve the matrix because this method is less expensive and almost as efficient as its variants with more parameters. As mentioned before, the HH run has not converged after 887 seconds but we used the same internal parameters in the SMT 2.0 software to compare the various models fairly. To be as exhaustive as possible, we also added the results for the EHH kernel, more details about the models without PLS are available in~\cite{Mixed_Paul}.

\subsubsection{Structural modelling: a cantilever beam bending problem \textcolor{black}{($n=2$, $ m=0$, $l=1$ and $L_1=12$) }}
\label{SMO_subsec:beam}

\begin{figure}[ht]
\centering
\vspace{-5pt}
\captionsetup{justification=raggedright,singlelinecheck=false}
\subfloat[Bending problem.]{
\begin{tikzpicture}
    \hspace{-3pt}
    %the points
    \point{origin}{-0.75}{-0.25};
    \point{begin}{0}{0};
    \point{end}{5}{0};
    \point{end_bot}{4.99}{-0.9};
    \point{end_up}{5}{0.5};
    %the beam
    \beam{2}{begin}{end};
    %the support
    \support{3}{begin}[-90];
    %the load
    \load{1}{end}[90]   ;
    %the inscription of the load
 %   \notation{1}{origin}{0};
    \notation{1}{end_up}{$F=50kN$};

     \draw[<->] (end) -- (end_bot) node[midway, right] {$\delta$} ;
     \draw[<->] (0,0.5) -- (5,0.5) node[midway, above] {L};
     
    %the deflection curves
    \draw
      [-, ultra thick] (begin) .. controls (1.5, +.01) and (2.5, -.15) .. (4.93, -0.9)
      [-, ultra thick] (begin) .. controls (1.5, +.01) and (2.5, -.2) .. (4.85, -1.5)
      [-, ultra thick] (begin) .. controls (1.5, +.01) and (2.5, -.4)   .. (4.78, -1.9);
  \end{tikzpicture}
\label{fig:beam}   
}
%\raggedright
%\hfill
\subfloat[Possible cross-section shapes.]{
\centering
\begin{tikzpicture}
\hspace{-50pt}

\vspace{-1.2cm}

%%%%%
\tstar{0.25}{0.5}{6}{0}{thick,fill=yellow,xshift=+3.6cm,yshift= -1.2cm}
\tstar{0.14}{0.28}{6}{0}{thick,fill=white,xshift=+3.6cm,yshift= -1.2cm}

\tstar{0.25}{0.5}{6}{0}{thick,fill=yellow,xshift=+2.4cm,yshift= -1.2cm}
\tstar{0.08}{0.16}{6}{0}{thick,fill=white,xshift=+2.4cm,yshift= -1.2cm}

\tstar{0.25}{0.5}{6}{0}{thick,fill=yellow,xshift=+1.2cm,yshift= -1.2cm}

%%%%%
\fill[green,even odd rule] (3.6,0) circle (0.5) (3.6,0) circle (0.33);
\draw (3.6,0) circle (0.5) ;
\draw (3.6,0) circle (0.33) 
; 
\fill[green,even odd rule] (2.4,0) circle (0.5)(2.4,0) circle (0.17);
\draw (2.4,0) circle (0.5) ;
\draw (2.4,0) circle (0.17) 
; 
\fill[green,even odd rule] (1.2,0) circle (0.5) ;
\draw (1.2,0) circle (0.5) ;

%%%%%
\def\pos{-2.4}
\fill[blue,even odd rule]  (\pos-0.5,-1.7+1.2) -- (\pos-0.5,-0.7+1.2) -- (\pos+0.5,-0.7+1.2) -- (\pos+0.5,-1.7+1.2) -- cycle ;

\def\pos{-1.2}
\fill[blue,even odd rule]  (\pos-0.5,-1.7+1.2) -- (\pos-0.5,-0.7+1.2) -- (\pos+0.5,-0.7+1.2) -- (\pos+0.5,-1.7+1.2) -- cycle   (\pos-0.25,-1.45+1.2) -- (\pos-0.25,-0.95+1.2) -- (\pos+0.25,-0.95+1.2) -- (\pos+0.25,-1.45+1.2) -- cycle ;

\def\pos{0}
\fill[blue,even odd rule]  (\pos-0.5,-1.7+1.2) -- (\pos-0.5,-0.7+1.2) -- (\pos+0.5,-0.7+1.2) -- (\pos+0.5,-1.7+1.2) -- cycle   (\pos-0.35,-1.55+1.2) -- (\pos-0.35,-0.85+1.2) -- (\pos+0.35,-0.85+1.2) -- (\pos+0.35,-1.55+1.2) -- cycle ;

%%%%%%
\def\pos{-2.4}
\fill[black] (\pos-0.24,-0.5-1.2) -- (\pos-0.24,0.5-1.2) -- (\pos+0.24,0.5-1.2)  -- (\pos+0.24,-0.5-1.2)   -- cycle ;
\fill[black] (\pos-0.5,-0.5-1.2) -- (\pos-0.5,-0.18-1.2) -- (\pos+0.5,-0.18-1.2)  -- (\pos+0.5,-0.5-1.2)   -- cycle ; 
\fill[black] (\pos-0.5,0.18-1.2) -- (\pos-0.5,0.5-1.2) -- (\pos+0.5,0.5-1.2)  -- (\pos+0.5,0.18-1.2)   -- cycle ;  

\def\pos{-1.2}
\fill[black] (\pos-0.19,-0.5-1.2) -- (\pos-0.19,0.5-1.2) -- (\pos+0.19,0.5-1.2)  -- (\pos+0.19,-0.5-1.2)   -- cycle ;
\fill[black] (\pos-0.5,-0.5-1.2) -- (\pos-0.5,-0.25-1.2) -- (\pos+0.5,-0.25-1.2)  -- (\pos+0.5,-0.5-1.2)   -- cycle ; 
\fill[black] (\pos-0.5,0.25-1.2) -- (\pos-0.5,0.5-1.2) -- (\pos+0.5,0.5-1.2)  -- (\pos+0.5,0.25-1.2)   -- cycle ;  

\def\pos{0}
\fill[black] (\pos-0.14,-0.5-1.2) -- (\pos-0.14,0.5-1.2) -- (\pos+0.14,0.5-1.2)  -- (\pos+0.14,-0.5-1.2)   -- cycle ;
\fill[black] (\pos-0.5,-0.5-1.2) -- (\pos-0.5,-0.32-1.2) -- (\pos+0.5,-0.32-1.2)  -- (\pos+0.5,-0.5-1.2)   -- cycle ; 
\fill[black] (\pos-0.5,0.32-1.2) -- (\pos-0.5,0.5-1.2) -- (\pos+0.5,0.5-1.2)  -- (\pos+0.5,0.32-1.2)   -- cycle ;  

\point{un}{-2.15}{-0.77};
\notation{1}{un}{\tiny 1};
\point{deux}{-2.15+1.2}{-0.77};
\notation{1}{deux}{\tiny 2};
\point{trois}{-2.15+2.4}{-0.77};
\notation{1}{trois}{\tiny 3};

\point{quatre}{-2.15+3.6}{-0.77};
\notation{1}{quatre}{\tiny 4};
\point{cinq}{-2.15+4.8}{-0.77};
\notation{1}{cinq}{\tiny 5};
\point{six}{-2.15+6}{-0.77};
\notation{1}{six}{\tiny 6};

\point{sept}{-2.15}{-0.77-1.2};
\notation{1}{sept}{\tiny 7};
\point{huit}{-2.15+1.2}{-0.77-1.2};
\notation{1}{huit}{\tiny 8};
\point{neuf}{-2.15+2.4}{-0.77-1.2};
\notation{1}{neuf}{\tiny 9};

\point{dix}{-2.15+3.6}{-0.77-1.2};
\notation{1}{dix}{\tiny 10};
\point{onze}{-2.15+4.8}{-0.77-1.2};
\notation{1}{onze}{\tiny 11};
\point{douze}{-2.15+6}{-0.77-1.2};
\notation{1}{douze}{\tiny 12};

  \end{tikzpicture}
  \label{fig:beam_shape}    
}
\captionsetup{justification=centering,singlelinecheck=false}
\caption{Cantilever beam problem~\cite[Figure 6]{Mixed_Paul}.}
\end{figure}
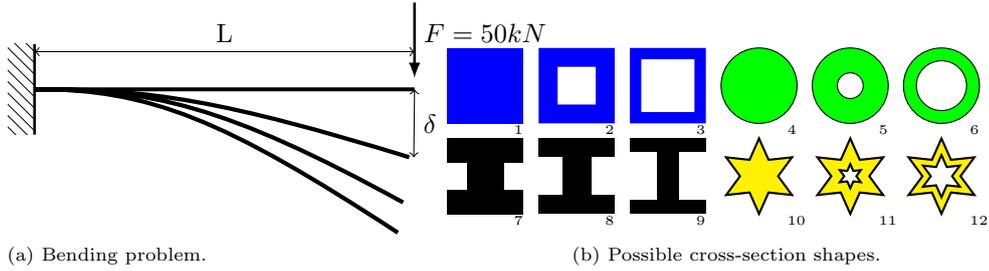

\begin{figure}[H]
\begin{center}

	\subfloat[GD kernel (3 hyperparameters: 1 cat. and 2 cont.)]{
      \centering 
		\includegraphics[clip=true, height=5.5cm, width=5cm]{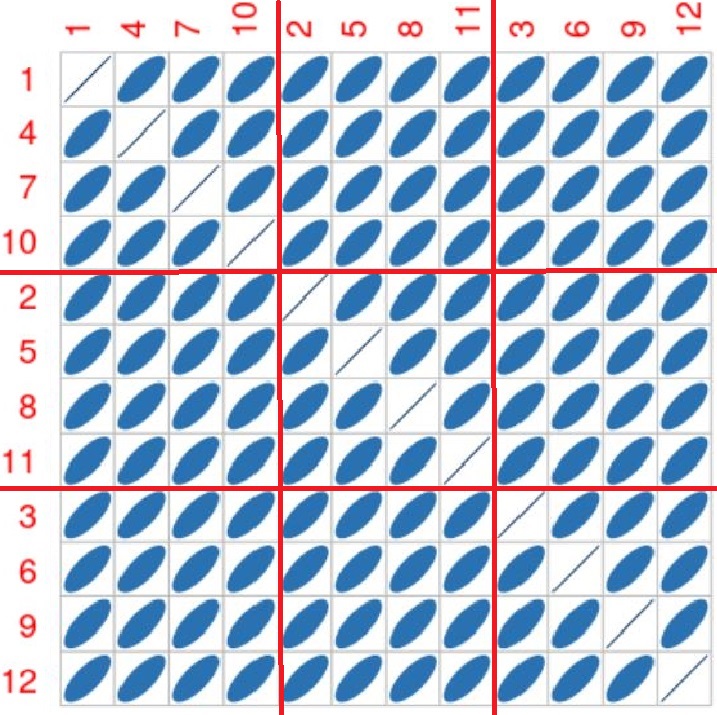} \label{SMO_corr_canti_gower}
     }  
              \hspace{.6 cm}
            \subfloat[HH with PLS kernel (3  hyperparameters: 1 cat. and 2 cont.)]{
      \centering
		\includegraphics[clip=true, height=5.5cm, width=5cm]{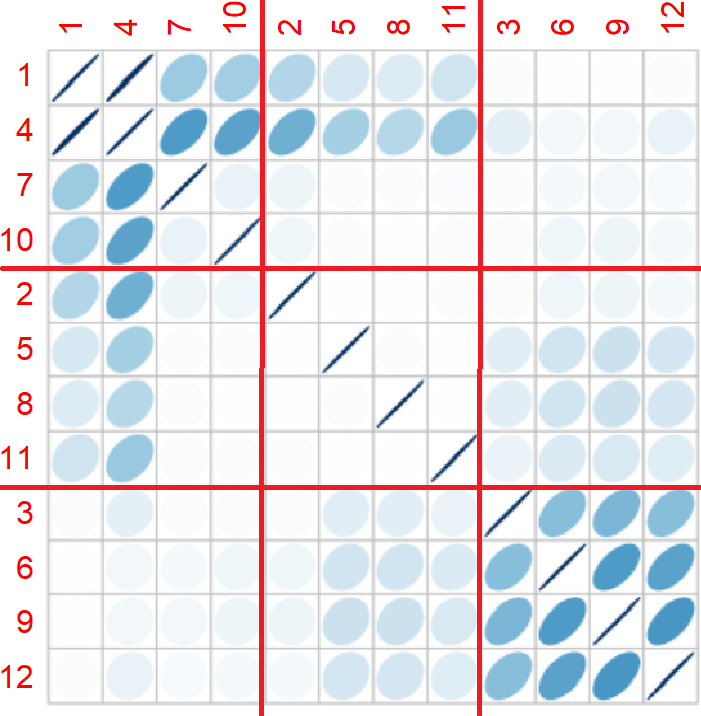}\label{SMO_corr_canti_hh_pls}
     }  
     \hspace{.1 cm}
      \centering
		\includegraphics[clip=true, height=5.5cm, width=1cm]{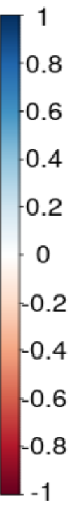}
     
        \subfloat[CR kernel (12 hyperparameters: 10 cat. and 2 cont.)]{
      \centering
		\includegraphics[clip=true, height=5.5cm, width=5cm]{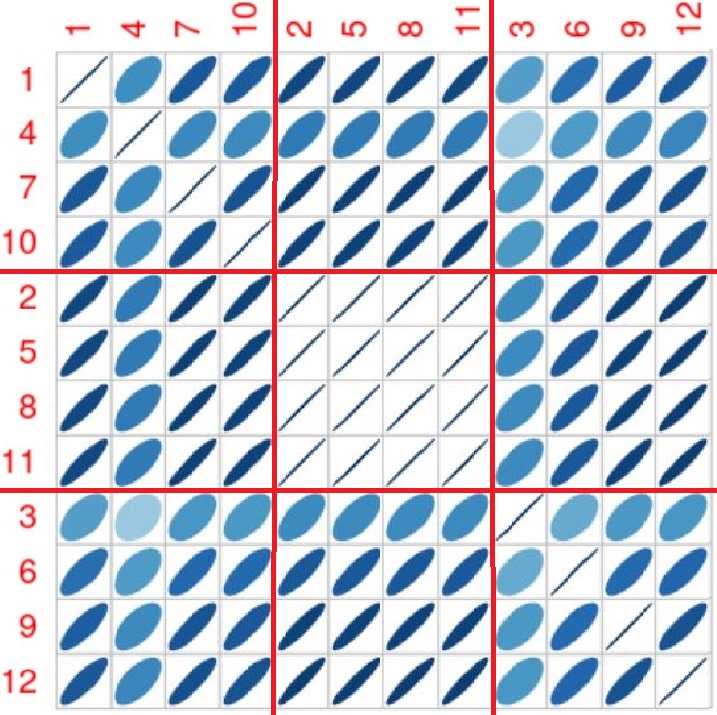} \label{SMO_corr_canti_cr}
     }  
        \hspace{.6 cm}
        \subfloat[HH kernel (66 hyperparameters: 64 cat. and 2 cont.)]{
      \centering
		\includegraphics[clip=true, height=5.5cm, width=5cm]{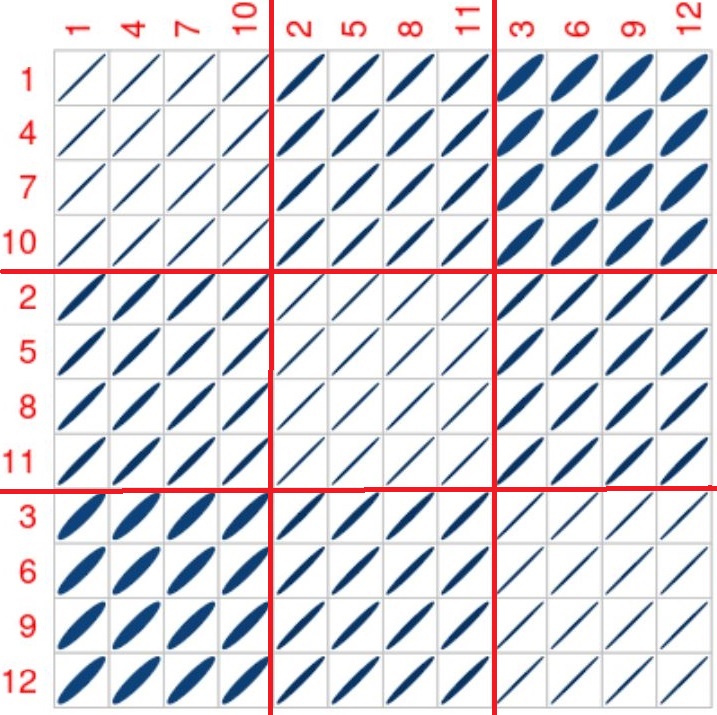}\label{SMO_corr_canti_ehh}
     }  
    \hspace{.1 cm}
      \centering
		\includegraphics[clip=true, height=5.5cm, width=1cm]{legend_pos_neg.png}
     
\caption{Correlation matrix $R_1^{cat}$  using different choices for $\Theta_1$ for the categorical variable $\tilde{I}$ from the cantilever beam problem.}
\label{SMO_corr_Cantilever}
\end{center} 
\end{figure}

A classic engineering problem frequently employed for model validation is the beam bending problem in its linear elasticity range~\cite{Roustant, Cheng2015TrustRB}. This problem serves as an illustrative example and involves a cantilever beam subjected to a load applied at its free end, denoted as $F$. The specific setup of the problem is depicted in~\figref{fig:beam}. In accordance with the findings presented in~\cite{ Cheng2015TrustRB}, the Young's modulus for the material is determined to be $E=200$ GPa, and a load of $F=50$ kN has been chosen for the analysis. Furthermore, following the methodology outlined in~\cite{Roustant}, a total of 12 potential cross-sections can be used for the beam. These cross-sections encompass four distinct shapes (square, circle, I and star), each with the possibility of being either full, thick, or hollow, as visually depicted in~\figref{fig:beam_shape}. For a given cross-section, which consists of a specific shape and thickness, its size is determined by the surface area denoted as $S$. Additionally, each cross-section is associated with a normalized moment of inertia $\tilde{I}$ around the neutral axis, representing a latent variable connected to the beam's shape~\cite{oune2021latent}. 
Hence, the problem at hand involves modeling with two continuous variables: the length $L$, ranging from 10 to 20 meters, and the surface area $S$, ranging from 1 to 2 square meters. Additionally, there is a categorical variable, $\tilde{I}$, with 12 levels representing the various cross-section options available. The tip deflection, at the free end, $\delta$ is given by $$ \delta = f( \tilde{I}, L,S) = \frac{F}{3E} \frac{L^3}{S^2\tilde{I}}. $$
\textcolor{black}{
As a result,  the relaxed dimension used to construct the GP model using the CR method is 14, while the relaxed dimension for the most general GP model employing the HH method is 68.} 
To compare our models, we draw a 98 point LHS as training set and the validation set is a grid of $12\times30\times30=10800$ points. For the four models GD, CR, HH and HH with PLS, the correlation matrix associated to every model are drawn in~\figref{SMO_corr_Cantilever} showing the predicted correlations between the available cross-sections. We recall that these matrices can be interpreted as such:
for two given levels $\{\ell_r^1,\ell_s^1\}$, the correlation term $[R_1]_{\ell_r^1,\ell_s^1}$ is in blue for correlation value close to 1, in white for correlations close to 0 and the thinner the ellipse, the higher the correlation.

The models are summarized in~\tabref{SMO_tab:resCantilever} indicating the complexity of each model and the information that could be recovered from it. For the HH kernel, the indicated computational time corresponds to the duration required to fully converge all 68 hyperparameters. In fact, the computational cost and difficulty to optimize the likelihood in spaces of dimension superior to ten is the biggest limitations of HH and such exhaustive kernels. This increase in difficulty to converge and associated computational cost is one of the main motivations for our PLS method and for simpler models.
As expected, we have 3 groups of 4 shapes depending on their respective thickness (respectively, the full levels \{1,4,7,10\}, the medium levels \{2,5,8,11\}, and the hollow levels \{3,6,9,12\}). The more the thickness is similar, the higher the correlation: the thickness has more impact than the shape of the cross-section on the tip deflection. However, given the database, two points with similar $L$ and $S$ values will have similar output whatever the cross-section. The effect of the cross-section on the output is always the same (in the form of $\frac{1}{\tilde{I}}$) leading to an high correlation after maximizing the likelihood. 
%
%For both squared exponential and absolute exponential kernels, the RMSE, likelihood and computational time for every model are shown in~\tabref{SMO_tab:resCantilever}. We recall that squared exponential and absolute exponential kernels differ only on the continuous variables and are the same for the categorical part.
%As expected, the computational time and the likelihood increase when the model is more complex. The DoE seems of sufficient size for this problem as the computed RMSE (\textit{i.e.}, the total displacement error) decreases with the model complexity.
%
%\resizebox{0.97\textwidth}{!}{%
\begin{table}[!ht]
\centering
\caption{Results of the cantilever beam models.}
\begin{tabular*}{\linewidth}{ccccc}
\hline
\textbf{Categorical kernel} &
%Displacement error (cm) &
\textcolor{black}{Identified clusters} & 
%L &
\# hyps.  &$\ $ Time (s) & \textcolor{black}{RMSE (cm)} \\
\hline 
 \textbf{GD}   &
 %1.3861 & 
 -- &
 %111.13&  
 3 %&{2.97} 
 & 8 
 & \textcolor{black}{1.3858}
\\  
\textbf{HH+PLS}   & 
%1.4674 &
\textbf{Hollow} cluster and \textbf{Full} cluster &
%95.35 & 
3 & 38
 & \textcolor{black}{1.2989}

\\
\textbf{CR}  &
%1.1671 &
\textbf{Medium} cluster & 
%155.32  & 
14  %& 41.02 \\
& 89
 & \textcolor{black}{1.1604}

\\
\textbf{HH}   &
%0.2033 &
\textbf{Full}, \textbf{Medium} and \textbf{Hollow} & 
%235.66 &
68  & 2128
 & \textcolor{black}{0.1247}

\\
\hline
\end{tabular*}
\label{SMO_tab:resCantilever}
\end{table}
%}
For the GD model in~\figref{SMO_corr_canti_gower}, there is only one mean positive correlation value, therefore no structural information can be extracted from this unique value. On the contrary, the HH model is the most general one and can model every cross-correlation value independently from the others and, in~\figref{SMO_corr_canti_ehh}, we can distinguish the three groups of four shapes, as expected, because the shapes of the cross-sections are not significant in comparison with their thickness.  Concerning the CR model, its structure favors the medium group, that is well represented, but the thick and hollow shapes are not distinguished in~\figref{SMO_corr_canti_cr}. 
To finish with, the HH+PLS model in~\figref{SMO_corr_canti_hh_pls} captures well the hollow shapes correlations and managed to capture half the full shapes correlations. More precisely, HH+PLS captured, with only one categorical hyperparameter (3 in total), the correlation between the full square and the full circle cross-sections (indexed 1 and 4) but failed to capture the correlation with full I and full star shapes (indexed 7 and 10). Consequently, with only one categorical hyperparameter, our model performs really well and is able to reconstruct structure from the data, thus outperforming the GD model for the same computational cost. 

To conclude, this section showed the capability of our PLS model to capture structures in the data while using only a small number of hyperparameters. Our model eases the optimization of the likelihood function and reduces the computational cost associated with the GP surrogate model. Notably, our work was efficiently applied to a structural modeling problem. However, building the GP surrogate is only part of the total optimization cost, and, in the next section, we show how our GP can be used in the context of high-dimensional MDO for aircraft design.

%Hence, CR can retrieve the three groups structure but can not obtain the high in-group correlations as the coefficients associated to both hollow and full group have to be smaller than the ones of the medium group. 

% To compare the mixed Kriging models, we draw a 98 point LHS as training set and the validation set is a grid of $12\times30\times30=10800$ points. 
% %
% For the four implemented methods, displacement error (computed with a root-mean-square error criterion), likelihood, number of hyperparameters and computational time for every model are shown in~\tabref{SMO_tab:resCantilever}.
% For the continuous variables, we use the square exponential kernel.
% More details are found in~\cite{Mixed_Paul}.
% As expected, the complex EHH and HH models lead to a lower displacement error and a higher likelihood value, but use more hyperparameters and increase the computational cost compared to GD and CR. 
% On this test case, the kernel EHH is easier to optimize than HH but in general, they are similar in terms of performance. 
% s\begin{table}[ht]
% \centering
% \caption{Results of the cantilever beam models.}
% \label{SMO_tab:resCantilever}
% \begin{tabular}{cccc}
% Categorical kernel &  Displacement error (cm) & $ \ $ Likelihood  & \# hyperparameters  %&$\ $ Time (s) \\
%   \\
%     \hline 
% GD   &1.3861 & 111.13&  3 %&{2.97} 
% \\  
%  CR   & 1.1671 & 155.32 & 14  %& 41.02 \\
% \\
% EHH  & {0.1613}   &
% {236.25} & 68 %& 1377 \\
% \\
% HH   & {0.2033} &
% {235.66} & 68 \\ % & 2128 \\
%   %HH+PLS & 4 & 27.284  & 21.95  \\
%   \hline 
% \end{tabular}
% \end{table}

 \subsection{Surrogate-based optimization: Bayesian optimization}
\label{SMO_sec:BO_val}

Efficient Global Optimization (EGO) is a well-known Bayesian optimizer that relies on GP to find out the optimum of an unconstrained black-box problem that can be evaluated a limited amount of times~\cite{Jones98}.
The workflow of EGO begins with building a first GP model based on an initial DoE, followed by employing an acquisition function to guide the selection of the next point that will be evaluated through the expensive black-box function. 
The most commonly used acquisition function is the expected improvement and, once a new point has been evaluated, the GP model is updated and the selection process repeats with the updated GP. At every step, a new model is built and a new point is evaluated until a maximal budget is reached. 
Hereinafter, we will use the GP models aforementioned to optimize expensive-to-evaluate black-box problems involving mixed integer variables using the EGO algorithm.
Moreover, EGO has been generalized to tackle constrained problems by Sasena \textit{et al.}~\cite{sasena2002exploration} with an algorithm called SEGO and used in the optimizer SEGOMOE~\cite{bartoli:hal-02149236}.

\subsubsection{Analytic validation on a mixed optimization problem \textcolor{black}{($n=1, \ m=0, \ l=1,$ and $L_1 = 10$)}} 
\label{SMO_sec:MI-BO}
\begin{figure}[H]
\begin{minipage}[b]{.6\linewidth}
\centering
\hspace{-1.25cm}
\subfloat[Convergence curves: medians of 20 runs.]{
\includegraphics[height=5cm,,width=7cm]
{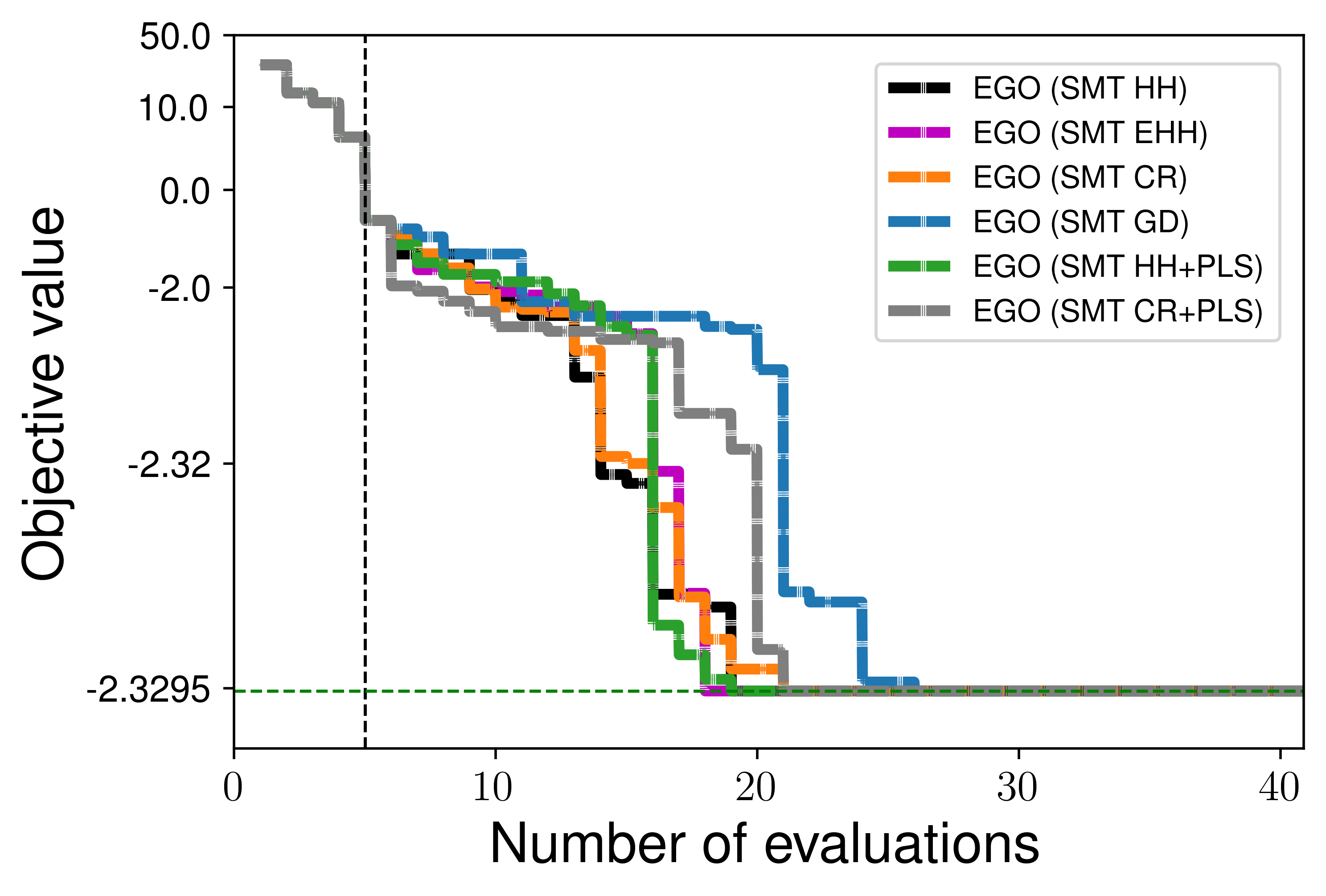}
\label{SMO_convmi}
}
\end{minipage}
\begin{minipage}[b]{.4\linewidth}
\centering 
\hspace{-1.5cm}
\subfloat[Boxplots after 25 evaluations.]{
\includegraphics[height=5cm,width=6.2cm]{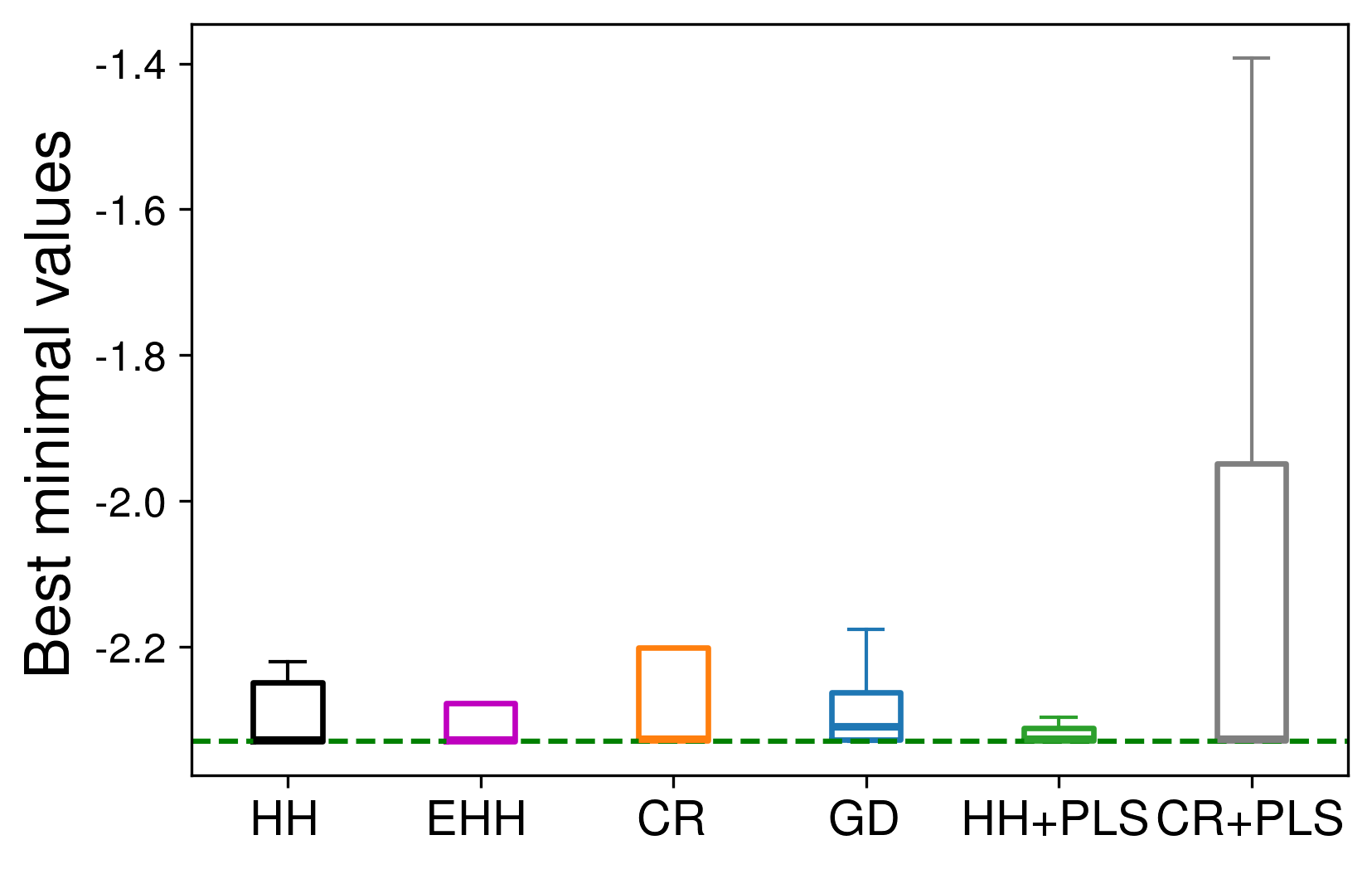}
\label{SMO_mini_mi}
}
\end{minipage}
\caption{Optimization results for the Toy function~\cite{CAT-EGO} for 20 DoE of 5 points.}
\label{SMO_res_optim_mi}
\end{figure}

\begin{figure}[H]
\begin{minipage}[b]{.6\linewidth}
\centering
\hspace{-1.25cm}
\subfloat[Convergence curves: medians of 20 runs.]{
\includegraphics[height=5cm,,width=7cm] {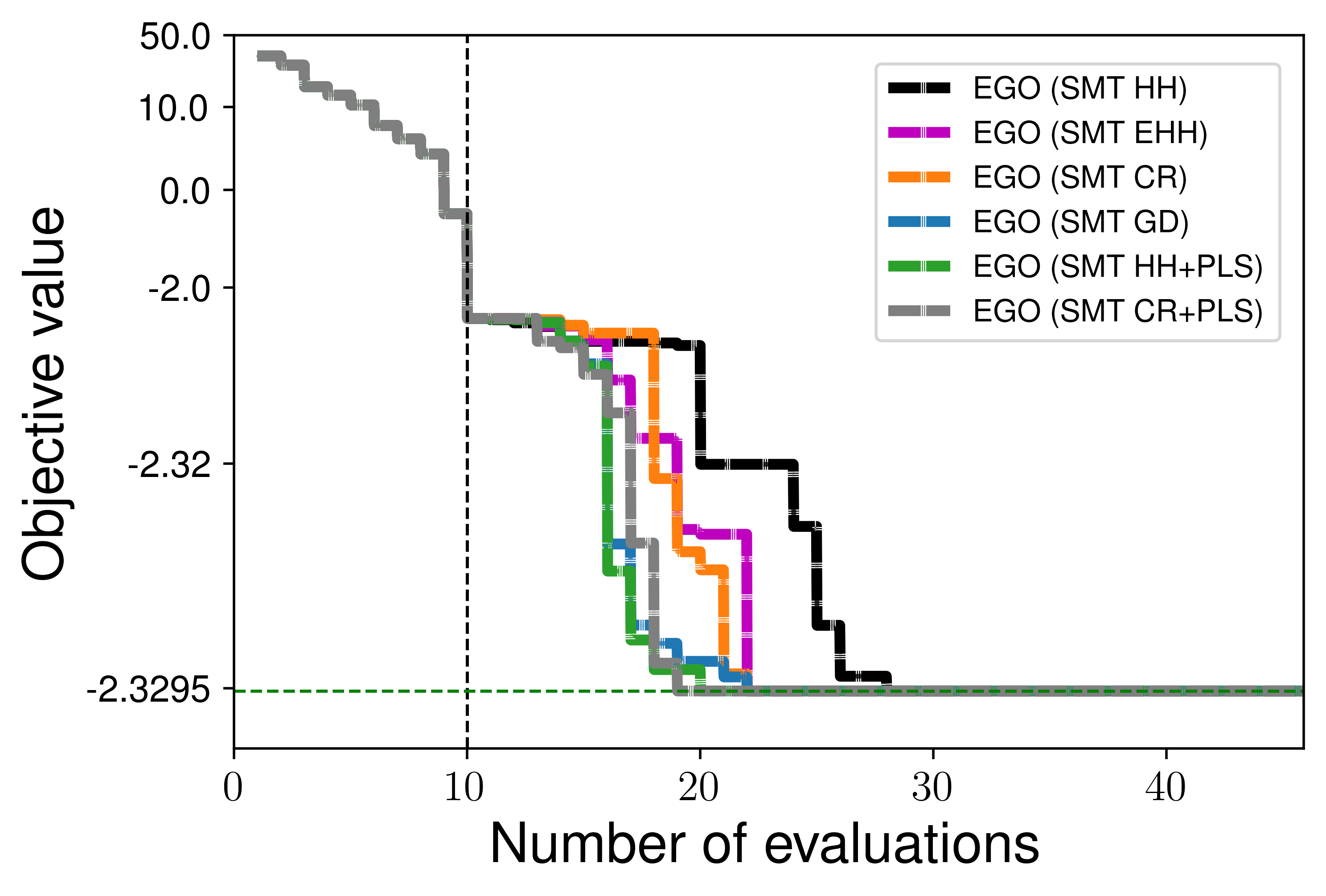}
\label{SMO_convmi2}
}
\end{minipage}
\begin{minipage}[b]{.4\linewidth}
\centering 
\hspace{-1.5cm}
\subfloat[Boxplots after 25 evaluations.]{
\includegraphics[height=5cm,width=6.2cm]{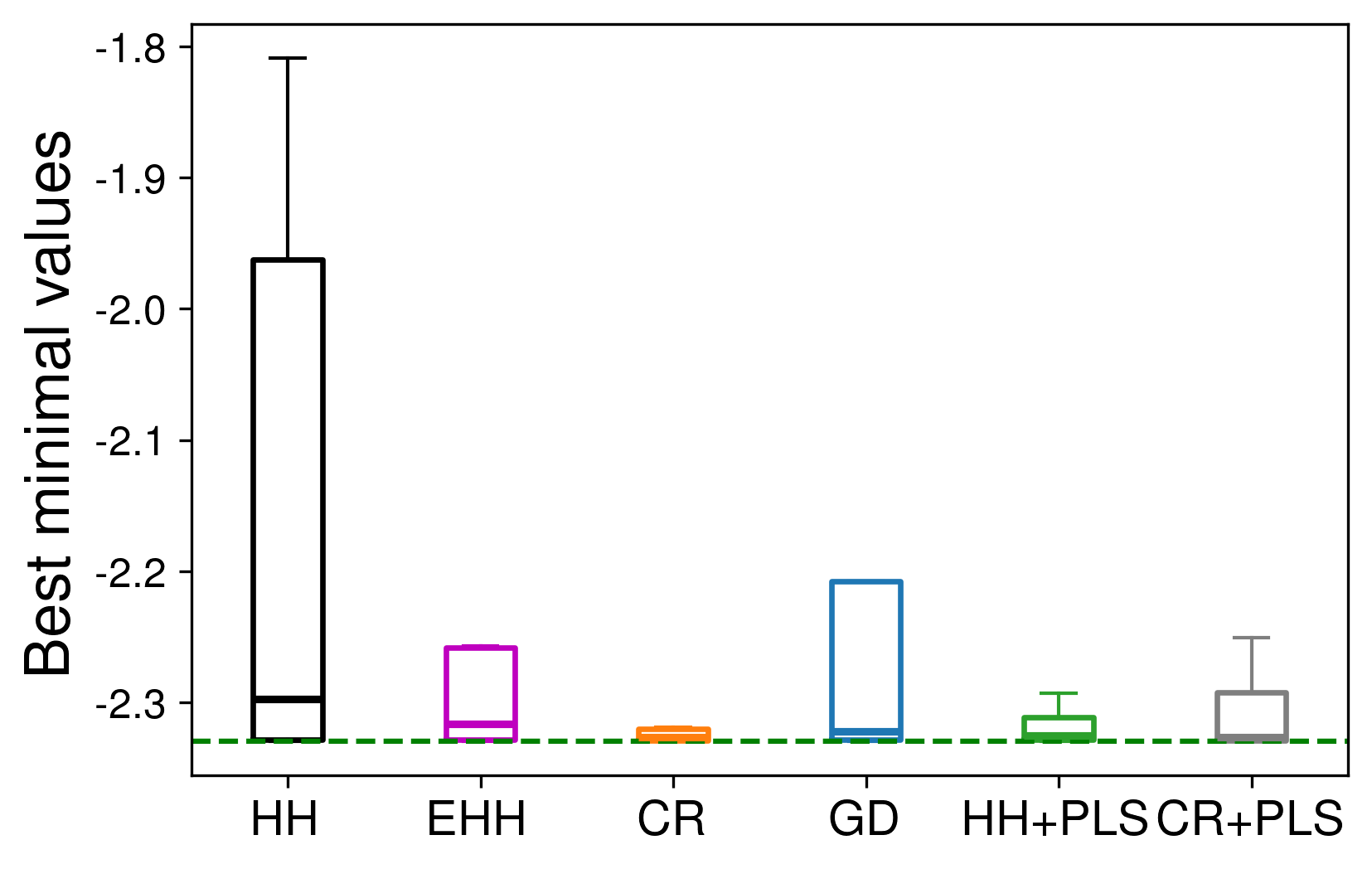}
\label{SMO_mini_mi2}
}
\end{minipage}
\caption{Optimization results for the Toy function~\cite{CAT-EGO} for 20 DoE of 10 points.}
\label{SMO_res_optim_mi2}
\end{figure}

The mixed test case that illustrates \textcolor{black}{BO} is a toy test case~\cite{CAT-EGO} detailed in Appendix~\ref{SMO_app:Toy}.
This test case has two variables, one continuous and one categorical with 10 levels.
\textcolor{black}{As a result, the relaxed dimension used to construct the GP model using the CR method is 11, while the relaxed dimension for the most general GP model employing the HH method is 46.} 
In~\figref{SMO_res_optim_mi} and~\figref{SMO_res_optim_mi2} six GP models are being compared. These six models are four classical 
\textcolor{black}{GP} %Kriging
models, namely GD, CR, EHH and HH and two PLS based models, namely our new model HH+PLS and the previously developed CR+PLS model~\cite{SciTech_cat}. 
In particular, in~\cite{SciTech_cat}, we used this CR+PLS method coupled with a criterion to choose automatically the number of PLS components that gives the best prediction for optimization. This adaptive PLS method for mixed integer has been applied to the MDO of DRAGON as detailed in Section~\ref{SMO_subsec:res_AD}.
To assess the performance of our algorithm, we performed 20 runs with different initial DoE sampled by LHS.
Every DoE consists of 5 points in~\figref{SMO_res_optim_mi} and of 10 points in~\figref{SMO_res_optim_mi2}. For both experiments, we chose a budget of 55 infill points. 
Figure \ref{SMO_convmi} and \figref{SMO_convmi2} plot the convergence curves for the six methods. To visualize the data dispersion, the boxplots of the 20 best solutions after 25 evaluations are plotted in~\figref{SMO_mini_mi} and~\figref{SMO_mini_mi2}. The computational times for every method are indicated in~\tabref{SMO_tab:resToy} for a 5 point DoE and in~\tabref{SMO_tab:resToy2} for a 10 point DoE. 
\textcolor{black}{We note that the overall computational cost is derived by the optimization cost related to the maximization of the infill criterion. In fact, such optimization is often related to the number of the design variables rather than the size of the DoE.}
On this test case, our method gives the best results with the 5 point DoE in terms of median convergence speed and dispersion among the 20 DoE. For the 10 point DoE, our method is among the faster together with CR+PLS. However, even if the HH+PLS method has been shown to be efficient for solving this test case, it is still more costly than CR+PLS or GD because the computational cost associated to the reconstruction of the matrix of hyperparameters is significant. 
Nevertheless, it is a method based only on two hyperparameters (one categorical and one continuous) making it around 20 times easier to optimize than HH or EHH and 3 times faster for better performance. In any case, using a 5 point DoE is slightly more efficient than using a 10 point DoE because \textcolor{black}{BO} is known to perform better with a smaller DoE for a given budget of evaluations~\cite{le2021revisiting}. But this effect is not significant for methods that use PLS, as PLS  benefits greatly from the initial DoE information to find the most interesting search directions. This explains why PLS methods are performing better with a 10 point DoE than with a 5 point DoE.
%As expected, the more a method is complex, the better the optimization. Both \texttt{SMT HH} and \texttt{SMT EHH} converged in around 18 evaluations whereas \texttt{SMT CR} and \texttt{SMT GD} take around 26 iterations as shown on~\figref{SMO_convmi}.
%Also, the more complex the model, the closer the optimum is to the real value as shown on~\figref{SMO_mini_mi}. 
%\vspace{-0.75cm}
\begin{table}[ht]    
\centering
 \caption{Results of the Toy problem optimization (5 point DoE and 55 infill points).}
\small
\small
\begin{tabular}{ccc}
\hline
\textbf{Kernel}& $ \ $ number of hyperparameters  & $\ $optimization duration (s) $\ $    \\
\hline 
\textbf{GD} & 2 &   315 \\
\hdashline 
\textbf{CR} & 11 &  503 \\
\textbf{CR+PLS} & 2 &  320 \\
\hdashline 
\textbf{HH}  & 46 &  1983 \\
\textbf{HH+PLS} & 2 &  646 \\
\hdashline 
\textbf{EHH} & 46 &  2086 \\
\hline
\end{tabular}
\label{SMO_tab:resToy}
\end{table}
%\vspace{-0.2cm}

%\vspace{-0.75cm}
\begin{table}[ht]    
\centering
 \caption{Results of the Toy problem optimization (10 point DoE and 55 infill points).}
\small
\small
\begin{tabular}{ccc}
\hline
\textbf{Kernel}& $ \ $ number of hyperparameters  & $\ $optimization duration (s) $\ $    \\
\hline 
\textbf{GD} & 2 &  314  \\
\hdashline 
\textbf{CR} & 11 & 479  \\
\textbf{CR+PLS} & 2 & 326  \\
\hdashline 
\textbf{HH}  & 46 &  2142 \\
\textbf{HH+PLS} & 2 & 662  \\
\hdashline 
\textbf{EHH} & 46 & 2079   \\
\hline
\end{tabular}
\label{SMO_tab:resToy2}
\end{table}

\subsubsection{Multidisciplinary design optimization for a green aircraft \textcolor{black}{
($n=10, \ m=0, \ l=2, L_1 = 17$ and $L_2 = 2$)}}
\label{SMO_subsec:res_AD}
\begin{figure}[H]
\begin{centering}
\includegraphics[height=5cm]{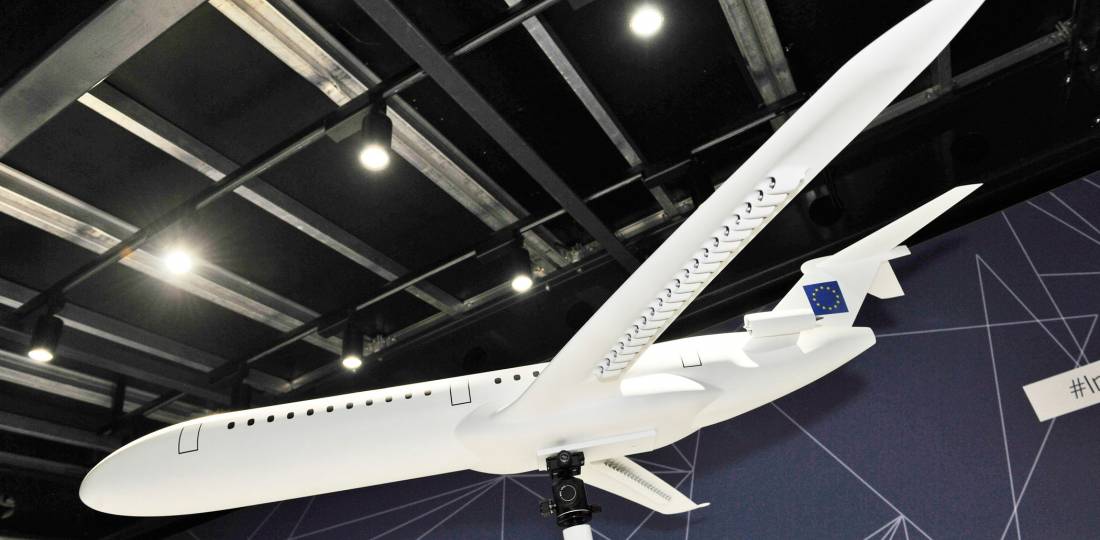}
   \caption{``\texttt{DRAGON}'' aircraft mock-up.}
        \label{SMO_Dragon2020}
\end{centering}
\end{figure}

For the core MDO application, we apply the Future Aircraft Sizing Tool with Overall Aircraft Design (FAST-OAD)~\cite{David_2021} on ``\texttt{DRAGON}'' 
(Distributed fans Research Aircraft with electric Generators by ONera), an innovative aircraft currently under development. The ``\texttt{DRAGON}'' aircraft concept in~\figref{SMO_Dragon2020} has been introduced by ONERA in 2019~\cite{schmollgruber} within the scope of the European CleanSky 2 program\footnote{\href{SMO_https://www.cleansky.eu/technology-evaluator}{\color{blue}https://www.cleansky.eu/technology-evaluator}} which sets the objective of 30\% reduction of CO2 emissions by 2035 with respect to 2014 state-of-the-art. A first publication in SciTech 2019~\cite{schmollgruber} was followed by an up-to-date report in SciTech 2020~\cite{schmollgruber2}.
In response to this ambitious goal, ONERA introduced a concept for a distributed electric propulsion aircraft that makes significant strides in enhancing fuel efficiency by optimizing propulsive performance. This is realized by an increase in the bypass ratio through a strategic placement of numerous compact electric fans on the wing pressure side, as an alternative to the use of larger turbofans. This design decision effectively resolves the challenges associated with large under-wing turbofans and grants the aircraft the capability to operate at transonic speeds. Consequently, the primary design objective for the  ``\texttt{DRAGON}'' revolves around accommodating a passenger capacity of 150 individuals and facilitating travel over a range of 2750 Nautical Miles at a speed of Mach 0.78.

\begin{figure}[H]
%\vspace*{-0.6cm}

  \centering 
\includegraphics[clip=true,,height=6cm]{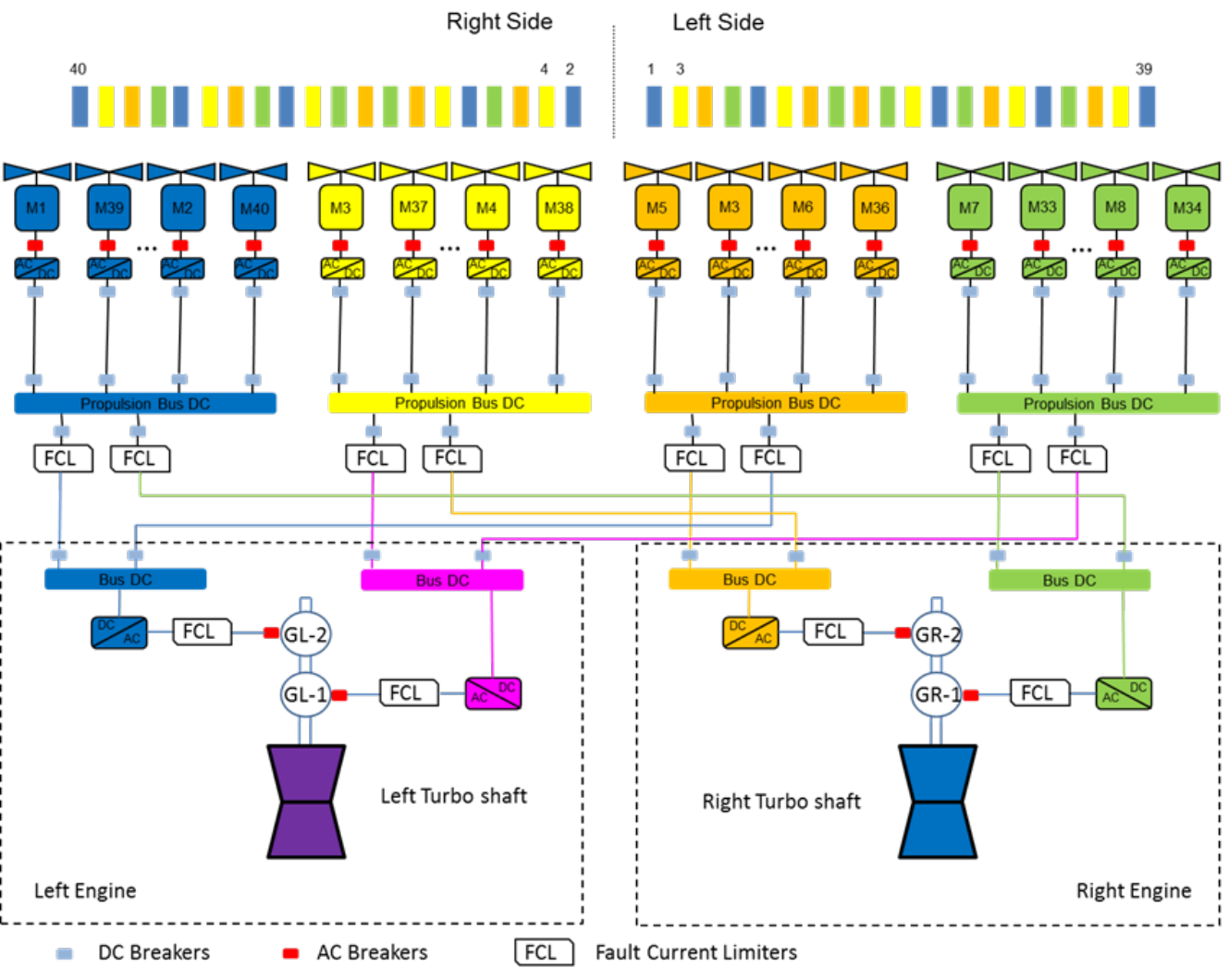}
 \caption{Turboelectric propulsive architecture.}
 \label{SMO_DragonArchitecture}
\label{SMO_Dragon}
\end{figure}

The integration of distributed propulsion in an aircraft introduces certain trade-offs. It necessitates the use of a turbo-electric propulsion system to provide the necessary power to drive the electric fans, which, in turn, contributes to increased intricacy and added weight. Typically, this power is generated onboard by coupling turboshafts to electric generators. The generated electrical power is subsequently transmitted to the electric fans through an electric architecture designed to ensure resilience in the face of potential single component failures. This safety feature is achieved through the deployment of redundant components, as illustrated  in~\figref{SMO_DragonArchitecture}.
The initial setup comprises two turboshafts, four generators, four propulsion buses with cross-feed capabilities, and 40 fans. This configuration was selected during the preliminary study phase due to its compliance with safety standards. Nevertheless, it was not specifically tailored for optimizing weight. Given that the turboelectric propulsion system significantly contributes to the overall weight of the aircraft, there is a specific interest in optimizing this system, especially concerning the number and type of individual components, each characterized by discrete or categorical values.
This time, as the evaluations are expensive, we are doing only 10 runs instead of 20.
Also, to have realistic results, the constraints violation will be forced to be less than $10^{-3}$. From now on, let MAC denote the Mean Average Chord, VT, the Vertical Tail, HT, the Horizontal Tail and TOFL, the Take-off Field Length. The optimizations are realized with SEGOMOE.

In~\cite{saves2021constrained,SciTech_cat}, the ``\texttt{DRAGON}'' configuration has been already optimized to attain such goals. In particular, in~\cite{SciTech_cat}, the model has been updated based on the results obtained in~\cite{saves2021constrained} that display limitations based on the turboshaft layout with the turbogenerators at the rear of the fuselage so it would be advantageous to locate the turbogenerator below the wing. 
Nonetheless, adopting this approach would impose constraints on the available space allocated for the electric fans. Consequently, it would restrict the maximum achievable propulsive efficiency. To address the inherent trade-off between a lighter propulsion system and an enhanced propulsive efficiency, we integrated the layout as a categorical variable within the optimization problem of~\figref{SMO_tab:dragon}. 
Finally, five constraints are considered on this optimization problem. Among these constraints, the TOFL, climb duration, and the top of climb slope angle exert a significant influence on the design of the hybrid electric propulsion system. Additionally, a portion of the wing trailing edge near the wingtip must be kept unobstructed to accommodate ailerons, thus limiting the available space for the electric fans. Lastly, compliance with wingspan limits is mandated by airport regulations.

To know how optimizing the fuel mass will impact the aircraft design, we are considering the optimization problem described in~\tabref{SMO_tab:dragon}. We can now solve a constrained optimization problem with 10 continuous design variables and 2 categorical variables with 17 and 2 levels respectively, for a total of 12 design variables. For the optimization, this new problem is a hard test case involving 29 relaxed variables and 5 constraints. The definition of the architecture variable is given in~\tabref{SMO_tab:dragon_archi1} and the definition of the turboshaft layout is given in~\tabref{SMO_tab:dragon_archi2}. 
\textcolor{black}{For modeling electric architectures, it is more efficient to represent the architectural choices using two integer variables instead of one categorical variable. 
However, taking this approach expands the range of potential architectures beyond the initial 17 configurations. Yet, there are two important constraints to consider for these possible setups.
The first constraint relates to the electrical connections between components: ensuring a certified electric architecture is crucial, and figuring out how to connect, for example, 8 motors to 6 generators is not straightforward.
The second constraint is connected to the distributed propulsion system, especially the numerous propellers. Managing this system involves addressing a substantial number of potential failures in the electro-mechanical architecture as for both stability and redundancy, not all electric connections are allowed.
Consequently, to simplify the optimization problem and avoid introducing many constraints, the model uses a single categorical variable to represent the various feasible architectures.}

Note that a simplified analysis has been done in a conference paper~\cite{SciTech_cat}, the latter was an optimization of the same aircraft configuration but with simpler methods, both HH and HH with PLS had never been tested before. 
\textcolor{black}{The relaxed dimension used to construct the GP model using the CR method is 29 as indicated in~\tabref{SMO_tab:dragon}, while the relaxed dimension for the most general GP model employing the HH method is 137.}

%\resizebox{0.97\textwidth}{!}{%
   % \hspace{-0.1cm}
\begin{table}[ht]
\centering
\vspace*{-0.3cm}
\caption{Definition of the ``\texttt{DRAGON}'' optimization problem.}
\begin{tabular}{lllrr}
& Function/variable & Nature & \# & Range\\
\hline
\hline
Minimize & Fuel mass & cont & 1 &\\
\hline
w.r.t & \mbox{Fan operating pressure ratio} & cont & 1 & $\left[1.05, 1.3\right]$ \\  
     & \mbox{Wing aspect ratio} & cont & 1 &    $\left[8, 12\right]$ \\
    & \mbox{Angle for swept wing} & cont & 1 & $\left[15, 40\right]$  ($^\circ$) \\
     & \mbox{Wing taper ratio} & cont & 1 &    $\left[0.2, 0.5\right]$ \\
     & \mbox{HT aspect ratio} & cont & 1 &    $\left[3, 6\right]$ \\
    & \mbox{Angle for swept HT} & cont & 1 & $\left[20, 40\right]$  ($^\circ$) \\
     & \mbox{HT taper ratio} & cont & 1 &    $\left[0.3, 0.5\right]$ \\
 & \mbox{TOFL for sizing}  & cont &1 & $\left[1800, 2500\right]$ ($m$) \\
 & \mbox{Top of climb vertical speed for sizing} & cont & 1 & $\left[300, 800\right]$ ($ft/min$) \\
 & \mbox{Start of climb slope angle} & cont & 1 & $\left[0.075, 0.15\right]$ ($rad$) \\
 & \multicolumn{2}{l}{Total  continuous variables} & 10 & \\
 \hline
& \mbox{Architecture} (levels) & cat & 17  & \{1,2, \ldots,16,17\} \\
& \mbox{Turboshaft layout}  (levels) & cat & 2  & \{1,2\} \\
 & \multicolumn{2}{l}{Total categorical variables} & 2 & \\
 \hline
  &   \multicolumn{2}{l}{\textbf{Total number of relaxed variables}} & {\textbf{29}} & \\
  \hline
subject to & Wing span \textless  \ 36   ($m$)  & cont & 1 \\
 & TOFL \textless  \ 2200 ($m$) & cont & 1 \\
 & Wing trailing edge occupied by fans  \textless  \ 14.4 ($m$) & cont & 1 \\
 & Climb duration \textless  \ 1740 ($s $) & cont & 1 \\
 & Top of climb slope \textgreater \ 0.0108 ($rad$) & cont & 1 \\
 & \multicolumn{2}{l}{\textbf{Total number of constraints}} & {\textbf{5}} & \\
\hline
\end{tabular}

\label{SMO_tab:dragon}
\end{table}
%}
\begin{table}[ht]
\centering
 \caption{Definition of the architecture variable and its 17 associated levels.}
\small
\begin{tabular}{ccc}
\hline
  \textbf{Architecture number} & number of motors & number of cores and generators\\
  \hline
  \textbf{1} & 8 &2 \\
  \textbf{2} & 12 &  2\\
  \textbf{3} & 16 &  2\\
  \textbf{4} &20 &2 \\
  \textbf{5} & 24 &  2\\
  \textbf{6} & 28 &  2\\
  \textbf{7} &32 & 2\\
  \textbf{8} & 36  & 2\\
  \textbf{9} & 40 &  2\\
  \textbf{10} & 8   & 4\\
  \textbf{11} & 16  & 4\\
  \textbf{12} & 24  & 4\\
  \textbf{13} & 32  & 4\\
  \textbf{14} & 40  & 4\\
  \textbf{15} & 12  & 6\\
  \textbf{16} & 24  & 6\\
  \textbf{17} & 36  & 6\\
\hline
\end{tabular}
\label{SMO_tab:dragon_archi1}
\end{table}
\begin{table}[!h]
\centering
%\vspace*{-0.3cm}
 \caption{Definition of the turboshaft layout variable and its 2 associated levels.}
\small
\begin{tabular}{cccccc}
\hline
  \textbf{Layout} & position & y ratio & tail & VT aspect ratio & VT taper ratio\\
  \hline 
  \textbf{1} & under wing &0.25 & without T-tail& 1.8 & 0.3 \\
  \textbf{2} & behind & 0.34 & with T-tail& 1.2 & 0.85\\
\hline
\end{tabular}
\label{SMO_tab:dragon_archi2}
\end{table}

\begin{figure}[H]
\begin{minipage}[b]{.6\linewidth}
\centering
\hspace{-1.25cm}
\subfloat[Convergence curves: medians of 10 runs.]{
\includegraphics[height=5cm,,width=7cm]
{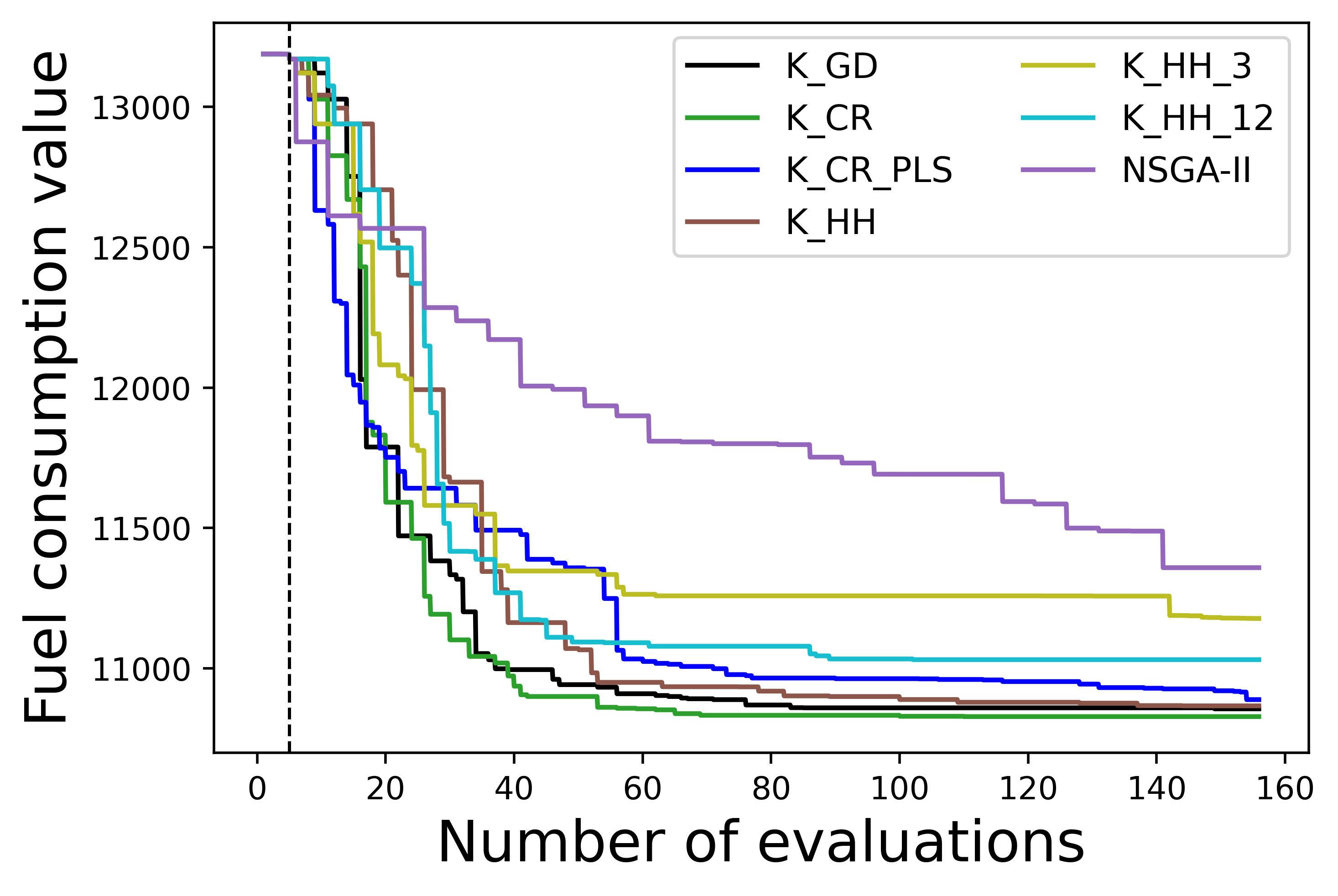}
\label{SMO_convmidragon5}
}
\end{minipage}
\begin{minipage}[b]{.4\linewidth}
\centering 
\hspace{-1.5cm}
\subfloat[Boxplots after 100 evaluations.]{
\includegraphics[height=5cm,width=6.2cm]{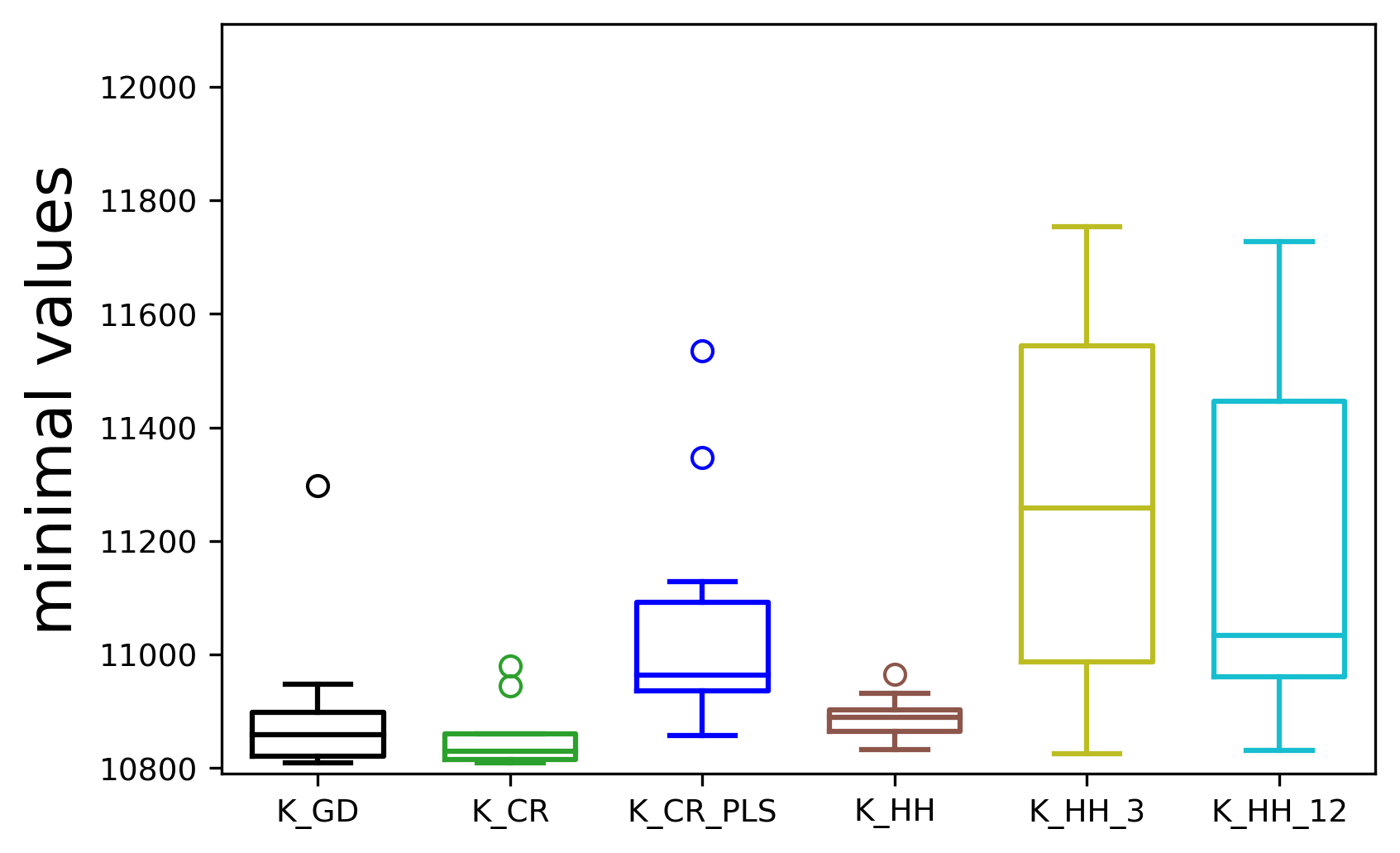}
\label{SMO_mini_midragon5}
}
\end{minipage}
\caption{Optimization results for the ``\texttt{DRAGON}'' aircraft~\cite{SciTech_cat}  for 10 DoE of 5 points.}
\label{SMO_res_optim_dragon5}
\end{figure}

\begin{figure}[H]
\begin{minipage}[b]{.6\linewidth}
\centering
\hspace{-1.25cm}
\subfloat[Convergence curves: medians of 10 runs.]{
\includegraphics[height=5cm,,width=7cm]
{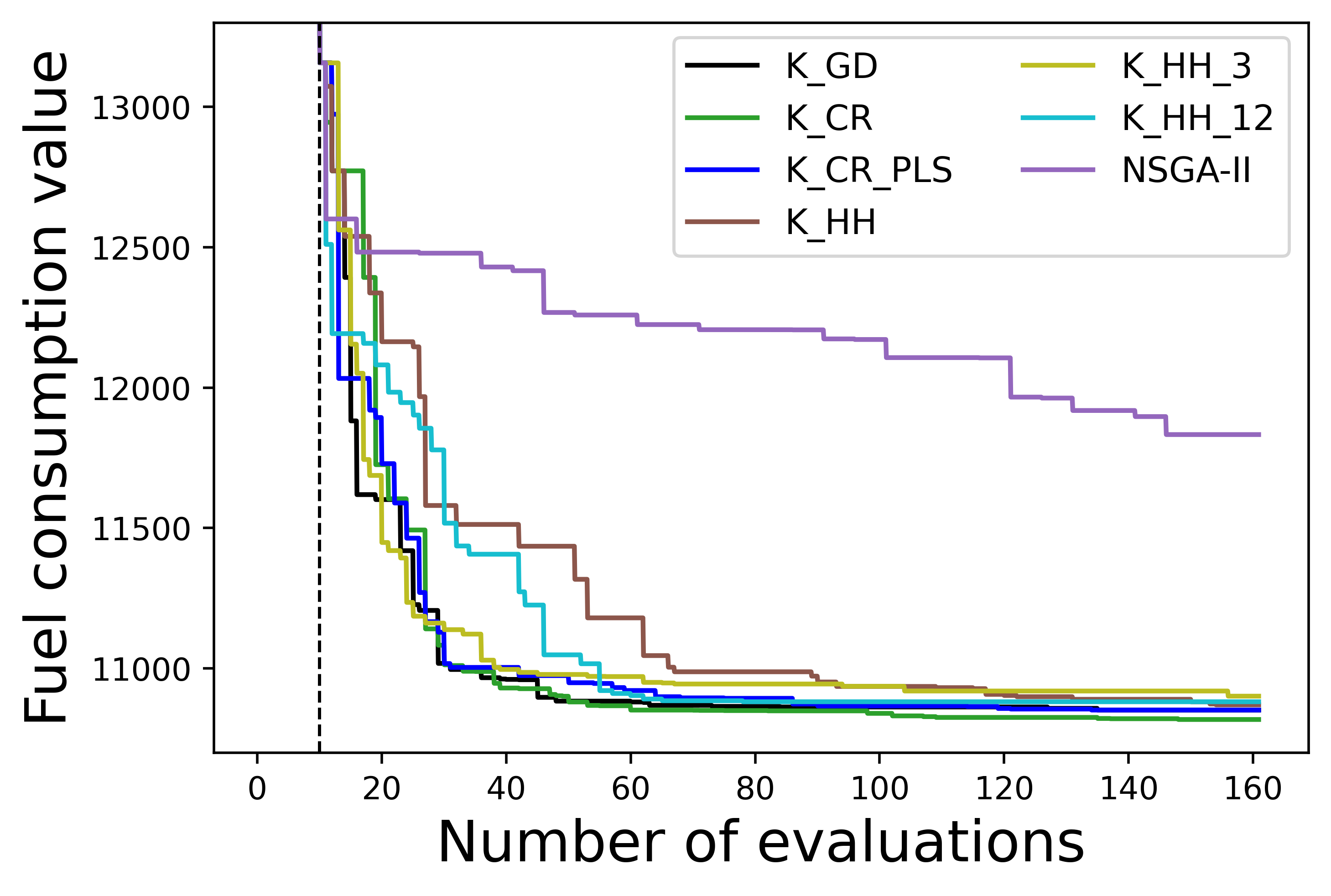}
\label{SMO_convmidragon10}
}
\end{minipage}
\begin{minipage}[b]{.4\linewidth}
\centering 
\hspace{-1.5cm}
\subfloat[Boxplots after 100 evaluations.]{
\includegraphics[height=5cm,width=6.2cm]{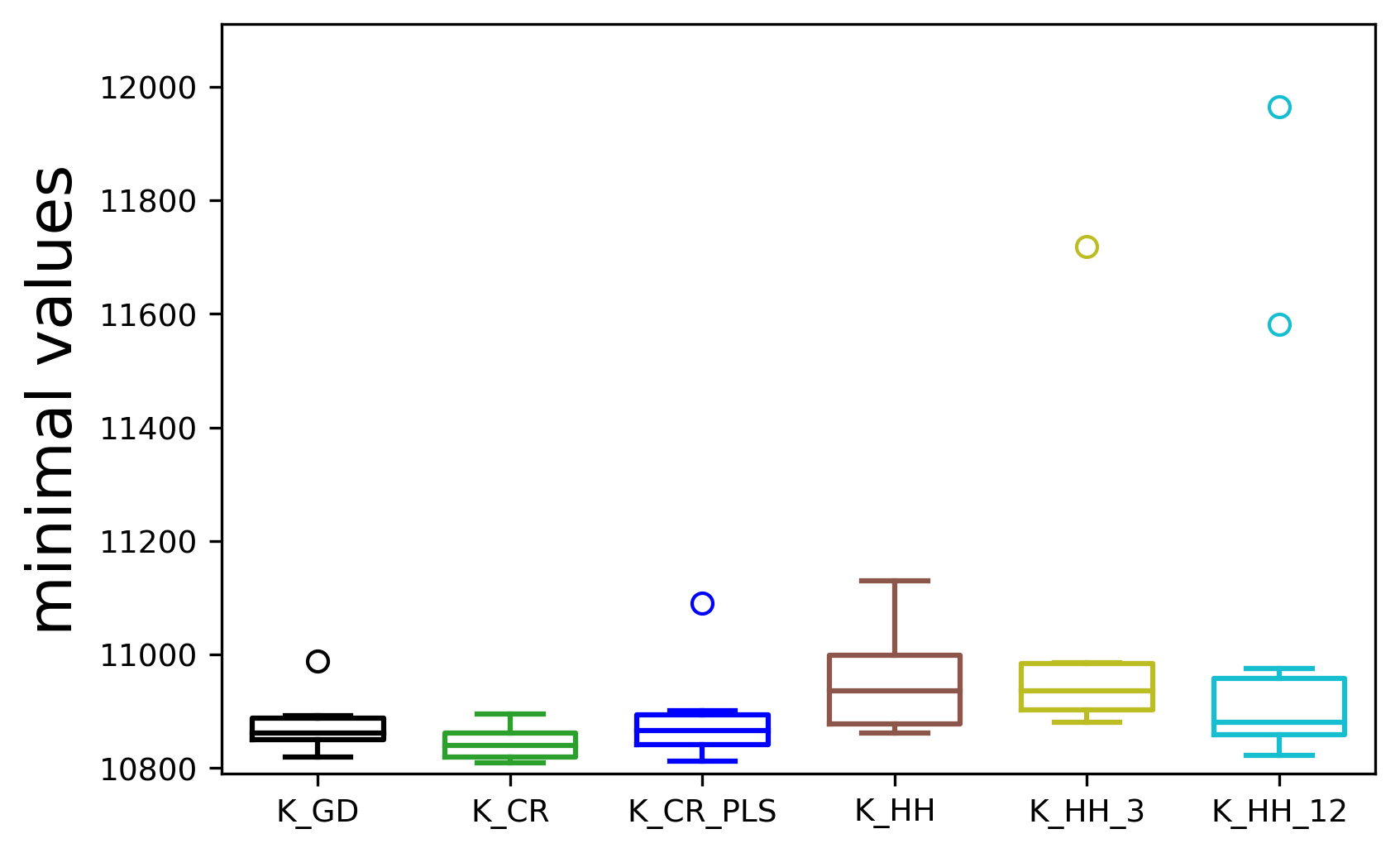}
\label{SMO_mini_midragon10}
}
\end{minipage}
\caption{Optimization results for the ``\texttt{DRAGON}'' aircraft~\cite{SciTech_cat}  for 10 DoE of 10 points.}
\label{SMO_res_optim_dragon10}
\end{figure}

To validate our method, we compare the 7 methods described in~\tabref{SMO_tab:dragon_meth} on the optimization of the ``\texttt{DRAGON}' aircraft concept with 5 and 10 points for the initial DoE as before.
\begin{table}[ht]
\centering
%\vspace*{-0.3cm}
 \caption{The various kernels compared on the MDO of ``\texttt{DRAGON}''. }
%\small
\begin{tabular}{cccc}
\hline
    \textbf{Name} & \# of cat. params & \# of cont. params & Total \# of params\\
    \hline 
    \textbf{GD} & 2 & 10 & 12\\
     \hdashline
    \textbf{CR} & 19 & 10 & 29  \\
    \textbf{CR with PLS 3D} & Not applicable & Not applicable & 3 \\ 
    \hdashline
    \textbf{HH} & 137 & 10 & 147 \\
    \textbf{HH with PLS 3D} & 2 & 1 & 3 \\
    \textbf{HH with PLS 12D} & 2 & 10 & 12 \\
  %  \textbf{EHH} & 137 & 10 & 147 \\
    \hdashline
    \textbf{NSGA-II} & Not applicable & Not applicable & Not applicable \\
    \hline
\end{tabular}
\label{SMO_tab:dragon_meth}
\end{table}   
As mentioned above, we are doing 10 runs for every method based on 10 starting DoE sampled by LHS to quantify the methods randomness. For every method and every starting DoE, we are running the method for a budget of 150 infill points, hence evaluating the black-box 155 times for the 5 point DoE and 160 times for the 10 point DoE. 
The results are given on~\figref{SMO_res_optim_dragon5} for the 5 point DoE and on~\figref{SMO_res_optim_dragon10} for the 10 point DoE.
More precisely,~\figref{SMO_convmidragon5} and~\figref{SMO_convmidragon10} display the convergence curves for the 7 methods and, to visualize the data dispersion, the boxplots of the 10 best solutions after 100 evaluations are shown in~\figref{SMO_mini_midragon5} and~\figref{SMO_mini_midragon10}. This computer experiments setup and figures are similar to what can be found in~\cite{SciTech_cat}.
\textcolor{black}{We note also that for this study, the computational cost of building the GP model is assumed to be negligible compared to the cost of evaluating the objective and the constraints at a given point. In fact, one simulation related to DRAGON is taking 2 to 5 minutes, which means half a day for a total of 150 simulations and around 60 days of computational time for running all the optimizations related to test case in this section.}

These results confirm the previous analyses made on the analytic test cases. 
First, the methods without PLS (GD, CR and HH) converge slightly better with a 5 point DoE, although this effect is less significant in this case because 5 or 10 points are both small quantities compared to 150 iterations and also because the search space is larger than before.
Second, the methods with PLS (CR-PLS, HH-PLS\_3D and HH-PLS\_12D) greatly benefit from a bigger starting DoE because the more representative the initial data, the more relevant the computed principal components. 
\textcolor{black}{
We note that a small DoE could lead to have PLS methods stuck in irrelevant zones as it might not be well-posed. In fact, in the case where an important part of the design space is not featured in the DoE, the available data will not be able to learn on this zone and thus the approximate representative space (computed by the PLS) will be sub-optimal.
}
The results show that for an initial DoE of 5 points, the methods without PLS are both faster and more consistent than their equivalent with PLS whereas the opposite is observed with the  DoE of 10 points. The 10 point DoE results also display that the HH model is costly and too complicated to be used efficiently with small DoE and that simpler methods like CR-PLS or GD are most effective at the start of the optimization process.
Based on this observation, it seems that combining models is the overall best method to tackle whatever black-box optimization problem beforehand.
These results also display that \textcolor{black}{BO} is more fitted to tackle such problems than evolutionary algorithms and that a small DoE could lead to have PLS methods stuck in irrelevant zones which, once again favors the alternative method of combining models for future works.

%%% A CHANGER%%
In terms of aircraft design, the ideal configuration was determined with an estimated fuel consumption of 10,809 kilograms against 11,248 kilograms for the original reference configuration. This configuration corresponds to option 10, featuring fewer engines (8 in total), but incorporating 4 cores and electric generators.
The most advantageous layout positions the turbo-generators at the rear. This choice is influenced by the increased lever arm between the wing and the horizontal tail, which results from the maximum sweep angle applied to the horizontal tail. However, it is worth noting that the combination of high sweep and high aspect ratio is not adequately considered from a structural standpoint, leading to unrealistically heavy weights for the horizontal stabilizer. Despite this limitation, the optimization process yields a suitable trade-off based on the models used in FAST-OAD.  The optimum found in the previous study of~\cite{SciTech_cat} was 10,816 kilograms which show that our new algorithms are more efficient. Still the new aircraft configuration is really close to the previous one, the changes are on the wing taper ratio reduced from 0.235 to 0.22, the TOFL for sizing which is now at the lower bound of 1800 m and not at 1803 m and to finish with, the start of climb angle was slightly reduced from 0.104 rad to  0.1035 rad. 

\section{Conclusion}
\label{SMO_sec:conclu}

% In this work, we have proposed a class of kernels for GP models that extends the exponential continuous kernels to the mixed-categorical setting. We showed that this class of kernels generalizes Gower distance and continuous relaxation based kernels. A classification between the proposed kernels as well as a proof of the SPD nature of the resulting correlation matrices have been also proposed. Numerical illustrations on analytical toy problems showed the good potential of the proposed kernels to reduce the number of hyperparameters and thus the computational time. The implementation of our proposed method has been released in the toolbox SMT v2.0\footnote{\url{https://smt.readthedocs.io/en/latest/}}. 

% When considering complex kernels, a good approach would be to use a model reduction technique such as Kriging with Partial Least Squares (KPLS)~\cite{Bouhlel18} that is derived from the construction of the correlation matrix via a kernel function. KPLS is an adaptation of the Partial Least Squares regression for exponential kernels and is used to reduce the number of hyperparameters and handle a large number of mixed inputs. Further works will consider to include such dimension reduction techniques to improve the computational efficiency of our model  and tackle higher dimensional problems. 

In this work, we proposed mixed-categorical metamodels based on GP for high-dimensional structural and multidisciplinary optimization. 
Our research was driven by the increasing complexity of engineering systems, which involve various disciplines and require optimization involving numerous design variables that could be either continuous, integer, and categorical. 
Our key findings center on the development of a more efficient approach for building surrogate models for large-scale mixed-categorical inputs. In~\cite{Mixed_Paul}, we introduced a mixed categorical kernel (EHH), a powerful tool for handling mixed-categorical variables which combine the matrix-based HH approach with the exponential kernel. However, we identified that the EHH and HH kernels effectiveness came at the cost of a significant increase in hyperparameters related to the GP surrogate model. To address this issue, we devised a novel approach by extending the partial least squares regression method as developed in~\cite{Bouhlel18} for continuous kernels, to reduce the number of hyperparameters while maintaining accuracy.

The significance of our research extends to both researchers and practitioners. For researchers, our work contributes to the evolving field of surrogate modeling for MDO. It offers a valuable solution to the challenge of high-dimensional mixed-categorical optimization, opening doors for further exploration in this domain.
Practitioners in engineering and optimization fields will find our findings beneficial as they provide a practical and efficient toolset for handling complex optimization problems. Our approach, implemented in the open-source software SMT~\cite{SMT2019,saves2023smt}, has been demonstrated effectively in structural and multidisciplinary applications, showcasing its real-world applicability.

%In conclusion, our study underscores the importance of addressing the challenges posed by mixed-categorical variables in high-dimensional multidisciplinary optimization. Through innovative methodologies and practical applications, we provided a valuable contribution to this field, enabling more efficient and accurate optimization processes for complex engineering systems.
% 
%Further works may include combining the several methods that now exist in the literature to have surrogate models that increase automatically in complexity when the size of the data set increases along the optimization process. 
%
\textcolor{black}{Further works may include combining the several methods that now exist in the literature to have surrogate models that increase automatically in complexity when the size of the dataset increases along the optimization process. 
Also, the surrogate models can be coupled to any surrogate-based optimization algorithm. In particular, in~\cite{grapin_constrained_2022}, SEGOMOE has been extended to multi-objective optimization and we also consider extending high-dimensional GP models to both mixed and hierarchical variables to tackle technological choices and variable-size problems~\cite{saves2023smt,audet2022general}.
Future performance benchmarks on industrial test cases should include comparisons with the Maximum Likelihood Estimation (MLE) approach for latent space identification, as demonstrated in latent map Gaussian process~\cite{oune2021latent}.
}

\section*{Acknowledgments}
This work is part of the activities of ONERA - ISAE - ENAC joint research group. The research presented in this paper has been performed in the framework of the COLOSSUS project (Collaborative System of Systems Exploration of Aviation Products, Services and Business Models) and has received funding from the European Union Horizon Programme under grant agreement n${^\circ}$ 101097120.
We thank Dr. Eric Nguyen Van (ONERA) and Christophe David (ONERA) for their contribution to DRAGON aircraft design.
We also thank Rémy Charayron for his constant help and support. 

\bmhead{Conflict of interest} 

The authors declare that they have no conflict of interest.

\bmhead{Replication of results} 
This section details the numerical tools used to obtain the results presented in this work. All numerical results were obtained using Python 3.8.8 and packages 
numpy 1.20.1, SciPy 1.6.2, pymoo 0.5.0, scikit-learn 1.0.2 and SMT: Surrogate Modelling Toolbox 2.0.
For constrained \textcolor{black}{BO}, we used the proprietary software SEGOMOE 5.0 by ONERA and ISAE-SUPAERO.
Default values of all functions were used unless stated otherwise.
%\textcolor{red}{Si SAMO il faudrait proposer les tutoriaux sur le site ...OBLIGATOIRE quand on prone l'opensource} 

\begin{appendices}

\section{Analytical modeling problem}
\label{SMO_subsec:cosine}
\begin{figure}[H]
\centering
\includegraphics[scale=.12]{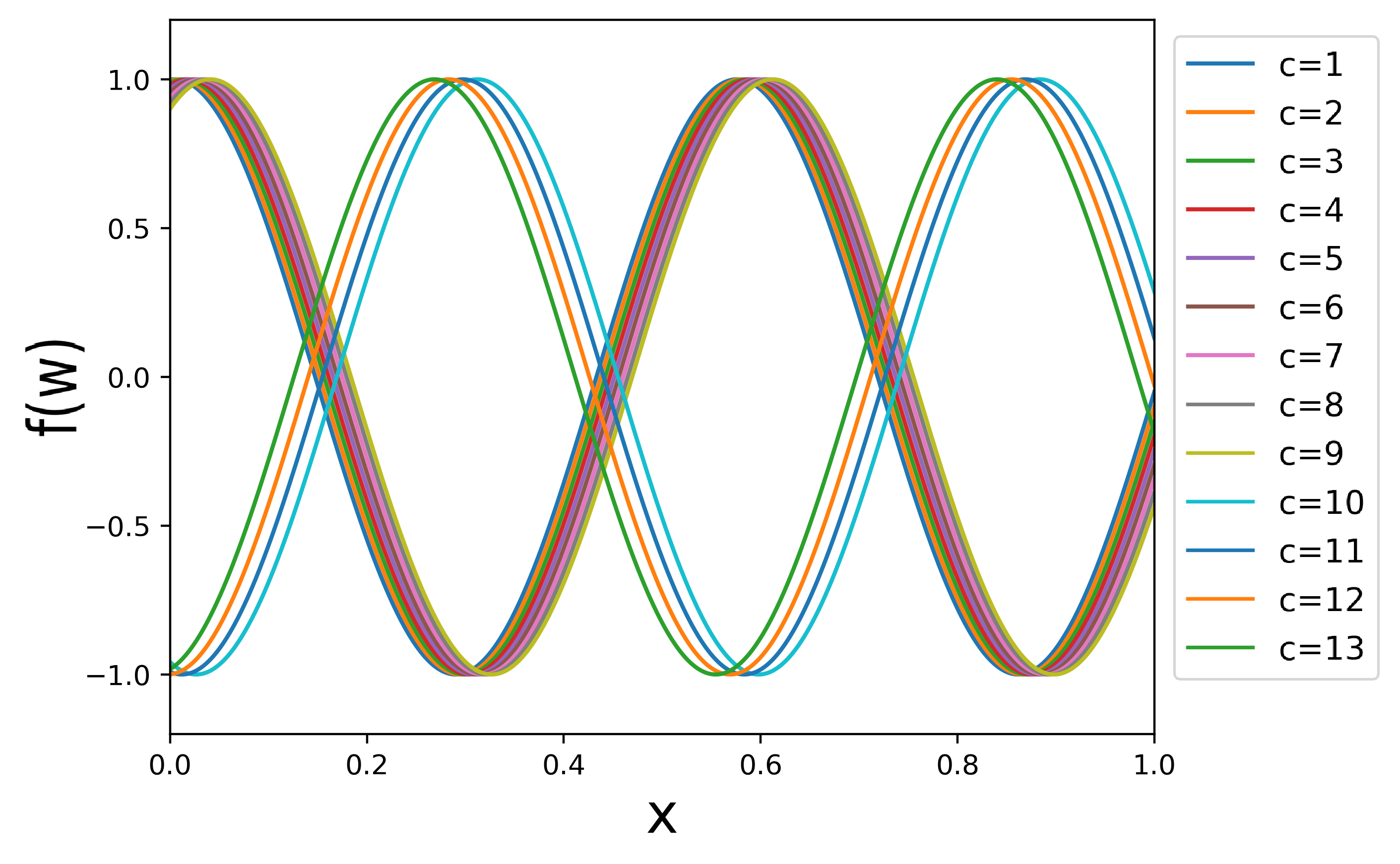}
\caption{Landscape of the cosine test case from~\cite{Roustant}.}
\label{SMO_fig:Roustant_ref}    
\end{figure}

This Appendix gives the detail of the categorical cosine test case of Section~\ref{SMO_subsec:catcos}.
This test case has one categorical variable with 13 levels and one continuous variable in $[0,1]$~\cite{Roustant}.
Let $w= (x,c )$ be a given point with  $x$ being the continuous variable and $c$ being the categorical variable, $c \in \{1, \ldots, 13\}$.

\begin{equation*}
\begin{split}
f(w) &= \cos \left( \frac{7 \pi}{2} x + \left( 0.4 \pi  + \frac{\pi }{15} c  \right) - \frac{c}{20} \right) , ~~~\mbox{if c $\in\{1,\ldots,9\}$ }  \\
f(w) &= \cos \left( \frac{7 \pi}{2} x  - \frac{c}{20} \right) , ~~~\mbox{if c $\in\{10,\ldots,13\}$ }  \\
\end{split}
\end{equation*}
The reference landscapes of the objective function (with respect to the categorical choices) are drawn on~\figref{SMO_fig:Roustant_ref}.

The DoE is given by a LHS of 98 points.
Our validation set is a evenly spaced grid of 1000 points in $x$ ranging  for every of the 13 categorical levels  for a total of 13000 points.

\section{Analytical optimization problem}
\label{SMO_app:Toy}
This Appendix gives the detail of the toy function of Section~\ref{SMO_sec:MI-BO}\footnote{\url{https://github.com/jbussemaker/SBArchOpt}}.
First, we recall the optimization problem:
\begin{equation*}
\begin{split}
& \min  f(w)=f(x, c_1) \\
& \mbox{w.r.t.} \ \  c_1 \in \{ 0,1,2,3,4,5,6,7,8,9 \} \\
& \quad \quad \quad x \in [ 0,1 ] \\
\end{split}
\end{equation*}
The objective function $f$ is defined as
\begin{equation*}
\begin{split}
   f({x}, {c_1})  =&   \mathds{1}_{c_1=0}  \ \cos(3.6 \pi(x-2)) +x -1 \\
     + &\mathds{1}_{c_1=1} \ 2 \cos(1.1 \pi \exp(x)) - \frac{x}{2} +2  \\
     + &\mathds{1}_{c_1=2}  \  \cos ( 2 \pi x)  + \frac{1}{2}x \\
  +&\mathds{1}_{c_1=3}  \  x ( \cos(3.4 \pi (x-1)) - \frac{x-1}{2})\\
  +& \mathds{1}_{c_1=4}  \   - \frac{x^2}{2} \\
 + &\mathds{1}_{c_1=5}  \  2 \cos(0.25 \pi \exp( -x^4))^2 - \frac{x}{2} +1 \\ 
 +& \mathds{1}_{c_1=6}   \ x \cos(3.4 \pi x ) - \frac{x}{2} +1 \\ 
 +&\mathds{1}_{c_1=7}   \  - x  (\cos(3.5 \pi x ) + \frac{x}{2}) +2 \\ 
 + &\mathds{1}_{c_1=8}   \ - \frac{x^5}{2} +1 \\ 
 +& \mathds{1}_{c_1=9}  \  - \cos (2.5 \pi x)^2 \sqrt{x} - 0.5 \ln (x+0.5)  - 1.3\\ 
\end{split}
\end{equation*}

\end{appendices}
\bibliography{sample.bib}
\end{document}